 
\overfullrule=0pt
\font\twelverm=cmr12
\font\twelvei=cmmi12
\font\twelvesy=cmsy10
\font\twelvebf=cmbx12
\font\twelvett=cmtt12
\font\twelveit=cmti12
\font\twelvesl=cmsl12
 
\font\ninerm=cmr9
\font\ninei=cmmi9
\font\ninesy=cmsy9
\font\ninebf=cmbx9
\font\ninett=cmtt9
\font\nineit=cmti9
\font\ninesl=cmsl9
 
\font\eightrm=cmr8
\font\eighti=cmmi8
\font\eightsy=cmsy8
\font\eightbf=cmbx8
\font\eighttt=cmtt8
\font\eightit=cmti8
\font\eightsl=cmsl8
 
\font\sixrm=cmr6
\font\sixi=cmmi6
\font\sixsy=cmsy6
\font\sixbf=cmbx6
 
\catcode`@=11 
\newskip\ttglue
 
\def\twelvepoint{\def\rm{\fam0\twelverm}
\textfont0=\twelverm  \scriptfont0=\ninerm
\scriptscriptfont0=\sevenrm
\textfont1=\twelvei  \scriptfont1=\ninei  \scriptscriptfont1=\seveni
\textfont2=\twelvesy  \scriptfont2=\ninesy
\scriptscriptfont2=\sevensy
\textfont3=\tenex  \scriptfont3=\tenex  \scriptscriptfont3=\tenex
\textfont\itfam=\twelveit  \def\it{\fam\itfam\twelveit}%
\textfont\slfam=\twelvesl  \def\sl{\fam\slfam\twelvesl}%
\textfont\ttfam=\twelvett  \def\tt{\fam\ttfam\twelvett}%
\textfont\bffam=\twelvebf  \scriptfont\bffam=\ninebf
\scriptscriptfont\bffam=\sevenbf  \def\bf{\fam\bffam\twelvebf}%
\tt  \ttglue=.5em plus.25em minus.15em
\normalbaselineskip=15pt
\setbox\strutbox=\hbox{\vrule height10pt depth5pt width0pt}%
\let\sc=\tenrm  \let\big=\twelvebig  \normalbaselines\rm}
 
\def\tenpoint{\def\rm{\fam0\tenrm}
\textfont0=\tenrm  \scriptfont0=\sevenrm  \scriptscriptfont0=\fiverm
\textfont1=\teni  \scriptfont1=\seveni  \scriptscriptfont1=\fivei
\textfont2=\tensy  \scriptfont2=\sevensy  \scriptscriptfont2=\fivesy
\textfont3=\tenex  \scriptfont3=\tenex  \scriptscriptfont3=\tenex
\textfont\itfam=\tenit  \def\it{\fam\itfam\tenit}%
\textfont\slfam=\tensl  \def\sl{\fam\slfam\tensl}%
\textfont\ttfam=\tentt  \def\tt{\fam\ttfam\tentt}%
\textfont\bffam=\tenbf  \scriptfont\bffam=\sevenbf
\scriptscriptfont\bffam=\fivebf  \def\bf{\fam\bffam\tenbf}%
\tt  \ttglue=.5em plus.25em minus.15em
\normalbaselineskip=12pt
\setbox\strutbox=\hbox{\vrule height8.5pt depth3.5pt width0pt}%
\let\sc=\eightrm  \let\big=\tenbig  \normalbaselines\rm}
 
\def\ninepoint{\def\rm{\fam0\ninerm}
\textfont0=\ninerm  \scriptfont0=\sixrm  \scriptscriptfont0=\fiverm
\textfont1=\ninei  \scriptfont1=\sixi  \scriptscriptfont1=\fivei
\textfont2=\ninesy  \scriptfont2=\sixsy  \scriptscriptfont2=\fivesy
\textfont3=\tenex  \scriptfont3=\tenex  \scriptscriptfont3=\tenex
\textfont\itfam=\nineit  \def\it{\fam\itfam\nineit}%
\textfont\slfam=\ninesl  \def\sl{\fam\slfam\ninesl}%
\textfont\ttfam=\ninett  \def\tt{\fam\ttfam\ninett}%
\textfont\bffam=\ninebf  \scriptfont\bffam=\sixbf
\scriptscriptfont\bffam=\fivebf  \def\bf{\fam\bffam\ninebf}%
\tt  \ttglue=.5em plus.25em minus.15em
\normalbaselineskip=11pt
\setbox\strutbox=\hbox{\vrule height8pt depth3pt width0pt}%
\let\sc=\sevenrm  \let\big=\ninebig  \normalbaselines\rm}
 
\def\eightpoint{\def\rm{\fam0\eightrm}
\textfont0=\eightrm  \scriptfont0=\sixrm  \scriptscriptfont0=\fiverm
\textfont1=\eighti  \scriptfont1=\sixi  \scriptscriptfont1=\fivei
\textfont2=\eightsy  \scriptfont2=\sixsy  \scriptscriptfont2=\fivesy
\textfont3=\tenex  \scriptfont3=\tenex  \scriptscriptfont3=\tenex
\textfont\itfam=\eightit  \def\it{\fam\itfam\eightit}%
\textfont\slfam=\eightsl  \def\sl{\fam\slfam\eightsl}%
\textfont\ttfam=\eighttt  \def\tt{\fam\ttfam\eighttt}%
\textfont\bffam=\eightbf  \scriptfont\bffam=\sixbf
\scriptscriptfont\bffam=\fivebf  \def\bf{\fam\bffam\eightbf}%
\tt  \ttglue=.5em plus.25em minus.15em
\normalbaselineskip=9pt
\setbox\strutbox=\hbox{\vrule height7pt depth2pt width0pt}%
\let\sc=\sixrm  \let\big=\eightbig  \normalbaselines\rm}
 
\def\twelvebig#1{{\hbox{$\textfont0=\twelverm\textfont2=\twelvesy
	\left#1\vbox to10pt{}\right.\n@space$}}}
\def\tenbig#1{{\hbox{$\left#1\vbox to8.5pt{}\right.\n@space$}}}
\def\ninebig#1{{\hbox{$\textfont0=\tenrm\textfont2=\tensy
	\left#1\vbox to7.25pt{}\right.\n@space$}}}
\def\eightbig#1{{\hbox{$\textfont0=\ninerm\textfont2=\ninesy
	\left#1\vbox to6.5pt{}\right.\n@space$}}}

 
\def\parag#1#2{\goodbreak\bigskip\bigskip\noindent
                   {\bf #1.\ \ #2}
                   \nobreak\bigskip\bigskip}
\def\intro#1{\goodbreak\bigskip\bigskip\bigskip\noindent
                   {\bf #1}\nobreak\bigskip\bigskip}
\long\def\th#1#2{\goodbreak\bigskip\noindent
                {\bf Theorem #1.\ \ \it #2}}
\long\def\lemma#1#2{\goodbreak\bigskip\noindent
                {\bf Lemma #1.\ \ \it #2}}
\long\def\prop#1#2{\goodbreak\bigskip\noindent
                  {\bf Proposition #1.\ \ \it #2}}
\long\def\cor#1#2{\goodbreak\bigskip\noindent
                {\bf Corollary #1.\ \ \it #2}}
\long\def\defin#1#2{\goodbreak\bigskip\noindent
                  {\bf Definition #1.\ \ \rm #2}}
\long\def\rem#1#2{\goodbreak\bigskip\noindent
                 {\bf Remark #1.\ \ \rm #2}}
\long\def\ex#1#2{\goodbreak\bigskip\noindent
                 {\bf Example #1.\ \ \rm #2}}
\def\proof{\nobreak\vskip.4cm\noindent{\it Proof.\ \ }}
\def\sqr#1#2{\vbox{\hrule height .#2pt
    \hbox{\vrule width .#2pt height #1pt \kern #1pt
      \vrule width .#2pt}
    \hrule height .#2pt }}
\def\square{\sqr74}
\def\endproof{\hphantom{MM}\hfill \llap{$\square$}\goodbreak}

\mathchardef\emptyset="001F
\mathchardef\hyphen="002D
\def\r{{\bf R}}

\def\Rm{\r^m}
\def\n{{\bf N}}

\def\q{{\bf Q}}
\def\vxi{\xi_1,\ldots,\xi_k}
\def\uxi{u_1,\ldots,u_k}
\def\cou{{\rm co}\{u_1(x),\ldots,u_k(x)\}}
\def\cox{{\rm co}\{\vxi\}}
\def\d{{\rm dom}}
\def\dF{\d F(\cdot,\om)}
\def\p{\prime}
\def\pp{\prime\prime}
\def\ris{{\rm ri}\Sigma_k}

\def\wh{\widehat}
 

\def\Extra{\noalign{\vskip3truept}}
 
\def\inf{\mathop{\rm\vphantom p inf}}
\def\Inf{\mathop{\rm inf}}
\def\Sup{\mathop{\rm sup}}
\def\s{\sum_{i=1}^k}
\def\min{\mathop{\rm min}}

\def\essup{\mathop{\rm ess\,sup}}

\def\mapright#1{\smash{\mathop{\longrightarrow}\limits_{#1}}}
 
\def\e{\varepsilon}
\def\f{\varphi}
\def\om{\Omega}
\def\g{\gamma}
 
\def\intb{\int_B}
\def\into{\int_\om}
 
\def\LL{L^{\infty}(\om ,\Rm)}
\def\LP{L^p(\om ,\Rm)}
\def\WP{W^{1,p}(\om ,\Rm)}
\def\WPP{W^{-1,p^\p}(\om ,\Rm)}
\def\W0{W^{1,p}_0(\om ,\Rm)}
\def\WPL{W^{1,p}(\om ,\Rm)\cap {L^\infty(\om ,\Rm)}}
\def\spc{C^1(\om)\cap W^{1,\infty}(\om)}
 
\def\F{\cal F}
\def\B{{\cal B}(\om)}
\def\A{{\cal A}(\om)}
\def\M{{\cal M}}
\def\Mu{{\cal M}_{\mu}}
 
\def\l{\lim_{h\rightarrow \infty}}
\def\li{\liminf_{h\rightarrow \infty}}
\def\Limsup{\mathop{\smash\limsup\vphantom\liminf}}
\def\lsr{\Limsup_{r\to\scriptstyle 0^+}}
 
\def\Eq{\,=\,}
\def\Ne{\,\ne\,}
\def\Ge{\,\ge\,}
\def\Le{\,\le\,}
\def\G{\,>\,}
\def\L{\,<\,}

\nopagenumbers
\magnification=1200
\hsize 14truecm
\vsize 25truecm
\hoffset=1.15truecm
\voffset=-1.45truecm
 
\null
\vskip 2truecm
{\twelvepoint
\baselineskip=1.7\baselineskip
\centerline{\bf INTEGRAL REPRESENTATION FOR A CLASS}
\centerline{\bf OF $\bf C^1$-CONVEX FUNCTIONALS}
}
\vskip2truecm
\centerline{Gianni DAL MASO}
\centerline{Anneliese DEFRANCESCHI}
\centerline{Enrico VITALI}
\vfil
 
{\eightpoint
\centerline{\bf Abstract}
\bigskip
\noindent
In view of the applications to the asymptotic analysis of a family of
obstacle problems, we consider a class of convex local functionals
$F(u,A)$, defined for all functions $u$ in a suitable vector valued
Sobolev space and for all open sets $A$ in ${\bf R}^n$. Sufficient
conditions are given in order to obtain an integral representation of the
form  $F(u,A)=\int_A f(x,u(x))\,d\mu + \nu(A)$, where $\mu$ and
$\nu$ are
Borel measures and $f$ is convex in the second variable.
 }
 
\vfil
\centerline {Ref. S.I.S.S.A.
71/M (May 92)}
\vskip 1truecm
\eject

\def\autori{\sl G.\ Dal Maso, A.\ Defranceschi, E.\ Vitali}
\def\titolo{\sl Integral representation for a class of
$C^1\!$-convex functionals}
\headline={\ifodd\pageno\rightheadline \else\leftheadline\fi}
\def\rightheadline{\eightpoint\hfil\titolo
\hfil\tenrm\folio}
\def\leftheadline{\tenrm\folio\hfil\eightpoint
\autori \hfil}
 
\pageno=1
\topskip=25pt
\vsize 22.5truecm
\hsize 16.2truecm
\hoffset=0truecm
\voffset=0.5truecm
\baselineskip=15pt
\hfuzz=2pt
\parindent=1cm
\mathsurround=1pt
 
\intro{Introduction}
 
This paper contains an integral representation theorem for a class of
convex local functionals which arise in the study of the asymptotic
behaviour of a sequence of minimum problems with obstacles for
vector valued Sobolev functions.
\par
Given an open subset $\om$ of $\r^n$ and $1<p<+\infty$, let $\WP$ be
the usual space of Sobolev functions with values in $\Rm$, and let
$\A$ be the family of all open subsets of~$\om$. The functionals
$F\colon \WP\times{\cal A}(\om)\rightarrow [0,+\infty]$ we are going
to consider are assumed to satisfy the following properties:
\smallskip
\item{(i)} (lower semicontinuity) for every $A\in{\cal A}(\om)$ the
function  $F(\cdot,A)$ is lower semicontinuous on~$\WP$;
\smallskip
\item{(ii)} (measure property) for every $u\in \WP$ the set function
$F(u,\cdot)$ is (the trace of) a Borel measure on~$\om$;
\smallskip
\item{(iii)} (locality property) $F(u,A)= F(v,A)$ whenever $u$,
$v\in\WP$, $A\in {\cal A}(\om)$, and $u|_A = v|_A$;
\smallskip
\item{(iv)} ($C^1\!$-convexity) for every $A\in {\cal A}(\om)$
the function $F(\cdot,A)$ is convex on $\WP$ and, in addition,
$
F(\f u + (1-\f)v,A)\Le F(u,A) + F(v,A)
$
for every $u$, $v\in\WP$ and for every $\f\in C^1(\om)\cap
W^{1,\infty}(\om)$ with $0\leq\f\leq 1$ on~$\om$.
\smallskip
This set of conditions is motivated by the study of the limit
behaviour, as ${h\to\infty}$, of a sequence of convex obstacle
problems of the form
$$
\min\{\into W(x,Du(x))\,dx : u\in H^1_0(\om,\Rm),\ u(x)\in K_h(x)\ \
\hbox{for q.e.\ }x\in A\}\,,\leqno{(0.1)}
$$
where $W(x,\xi)$ is quadratic with respect to~$\xi$, $A$ is an open
subset of $\om$ with $A\subset\subset\om$, and $K_h(x)$ is a closed
convex subset of $\Rm$ for every $h\in\n$ and for every $x\in\om$. By
using
\hbox{$\Gamma$-con}\-vergence techniques it is possible to prove
(see~[17]) that the limit problem can always be written in the
form
$$
\min\{\into W(x,Du(x))\,dx\,+\,F(u,A) : u\in H^1_0(\om,\Rm)\}\,,
$$
with $F$ satisfying the conditions considered above.
\par
In this paper we are concerned only with the properties of $F$ that
can be deduced {}from (i)--(iv). The main result (Theorem~6.5) is that
every functional $F$ satisfying (i)--(iv) can be written in the form
$$
F(u,A)\Eq \int_A f(x, u(x))\,d\mu \, + \, \nu(A)\,,\leqno(0.2)
$$
where $\mu$ and $\nu$ are positive Borel measures and
$f\colon\om\times\Rm\to[0,+\infty]$ is a Borel function, convex and
lower  semicontinuous in the second variable.
\par
This result will be used in a forthcoming paper (see~[17]) to
provide a detailed description of the structure of the limits of
sequences of obstacle problems of the form (0.1) under various
assumptions on $W(x,\xi)$ and~$K_h(x)$.
\par
Conditions (i)--(iii) are not enough to obtain an integral
representation of the form (0.2). Indeed, even convex functionals
depending on the gradient, like $F(u,A)=\int_A |Du|^pdx$, satisfy
(i)--(iii). Condition (iv) is  the most important one, and is
responsible for an integral  representation of the form (0.2), i.e.,
without terms depending on the gradient. This notion  of convexity,
also used, e.g., in [27] and [16], is strictly  related to the
notion of $C^1\!$-stability introduced by G. Bouchitt\'e and M. Valadier
in~[7], whose results are frequently used in our paper.
\par
For a general survey on integral representation theorems in $L^p$,
$W^{1,p}$, and $BV$ we refer to~[9]. See also [1], [2],
[6], [4], [21], [3] for more recent results.
\par
In the scalar case (i.e., $m=1$), integral representations on
$W^{1,p}(\om)$ of the
form (0.2), connected with limits of obstacle problems, can be found in
[18], [15], [5], [12], [16] under suitable
convexity conditions, and in [13] under monotonicity assumptions.
\par
Although the final statement in the vector case is exactly the same
as in the  scalar case, the proof is completely
different, since all arguments used in the papers mentioned above rely
on the order structure of~$\r$, involving truncations and monotonicity
methods.
\par
The main tools for the proof in the vector case are some technical
results obtained in our previous paper [14], based on the
theory of Lipschitz parametrization of convex sets  developed in [24]
and~[26]. In particular we shall use the following result (Theorem~2.9
in~[14]): given a finite number of  functions $u_1,\ldots, u_k$
in  $\WP\cap\LL$, their convex combinations with smooth  coefficients
form a dense subset in the set of all $W^{1,p}$-selections of the
polyhedral multivalued function $x\mapsto\cou$, where $\rm co$
denotes
the convex hull.
\par
The first step (Theorem~3.7) of our
result deals with the integral representation of the functional $F$ on
the set of all $W^{1,p}$-selections of such polyhedral multifunctions.
\par
In Theorem~5.4 we extend the integral representation of  $F$
to all the functions of $\WPL$ which satisfy, up to sets of capacity
zero, a suitable ``obstacle condition'' of the form $u(x)\in K(x)$,
which is necessary (but not sufficient) for the finiteness of the
functional. We note, incidentally, that the main difficulty in the
proof of our result lies in the fact that the functional $F$ is not
assumed to be finite everywhere, in view of the applications to
obstacle problems.
\par
The restriction to $\LL$, originated by the need of taking
products of $W^{1,p}$-functions, is dropped in Section~6. Moreover,
the ``obstacle  condition", given up to sets of capacity zero, is shown
to be equivalent to the condition $u(x)\in K(x)$ almost everywhere with
respect to a suitable measure (Proposition~6.3). This allows us to
obtain the integral representation (0.2) for every $u\in\WP$ and for
every $A\in\A$.
\par
In the last section (Theorem~7.3) we
prove that, if $F$ is quadratic or
positively $p$-homogeneous, then so is~$f$.
 
\intro{Acknowledgements}
This work is part of the Project EURHomogenization --
ERB4002PL910092
of the Program SCIENCE of the Commission of the European
Communities,
and of the Research Project ``Irregular Variational Problems'' of the
Italian National Research Council. The second and third author were
partially supported by a National Research Project of the Italian
Ministry of Scientific Research.
\bigskip
 
\parag{1}{Notation and preliminaries}
 
Throughout this paper
$m$, $n$ are two fixed positive integers, $p$ is a fixed real number,
$1<p<+\infty$, and $\om$ is an open subset of $\r^n$, possibly
unbounded. We
shall denote by $\A$ the family of the open subsets of $\om$ and by
$\B$ the
family of its Borel subsets. If $B\subseteq \r^n$ is a Borel set we denote
its
Lebesgue measure by $|B|$. The notation $a.e.$ stands for almost
everywhere with
respect to the Lebesgue measure.
\par
If $d$ is a positive integer, for every $x\in\r^d$ and $r>0$ we set
$B_r(x)= \{y\in\r^d:|y-x|<r\}\,,$
while ${\overline B}_r(x)$ denotes the closure of $B_r(x)$.
The $(d-1)$-dimensional simplex $\Sigma_d$
is defined by
$$
\Sigma_d\Eq\{\lambda\in \r^d: {\lambda^1+\cdots +\lambda^d\Eq
1,} \lambda^i\Ge 0\}\,,
$$
where $\lambda = (\lambda^1,\ldots,\lambda^d)$. If $C$ is a convex
subset of
$\r^d$, we denote by ${\rm ri\,}C$ its relative interior and by $\partial
C$
its relative boundary. In particular,  ${\rm ri\,}\Sigma_d=
\{\lambda\in\r^d:\lambda^1+\cdots+\lambda^d=1, \lambda^i>0\}$ and
$\partial\Sigma_d=\Sigma_d\setminus{\rm ri\,}\Sigma_d$.
\par
The
space $L^p(\om,\Rm)$ is endowed with the usual norm
$$
\|u\|_{L^p(\om,\Rm)}\Eq (\into|u|^pdx)^{1/p}\,.
$$
Let $\WP$ be the Banach space of all the functions $u\in L^p(\om,\Rm)$
with first
order distributional derivative $Du$ in $L^p(\om,\r^{mn})$, endowed
with the
norm
$$
\|u\|_{\WP}\Eq ({\| u\|}^p_{L^p(\om,\Rm)} +
{\|Du\|}^p_{L^p(\om,\r^{mn})})^{1/p}\,.
$$
The closure of $C^1_0(\om,\Rm)$ in $\WP$ will be denoted by $\W0$
($\Rm$ will
be omitted if $m=1$).
\medskip
\noindent
{\it Capacity.} For every compact set $K\subseteq\om$ we define the
capacity of
$K$ with respect to $\om$ by
$$
{\rm cap}(K,\om)\Eq {\rm
inf}\{\|\f\|^p_{W^{1,p}(\om)}:\f\in C^\infty_0(\om)\,, \f\Ge 1\
\hbox{on}\
K\}\,.
$$
The definition is extended to all subsets of $\om$ as external capacity
in the usual way (see, for example, [11] and [31]).
\par
Let $E$ be a subset of $\r^n$. If a statement depending on $x\in\r^n$
holds
for every $x\in E$ except for a set $N\subseteq E$ with capacity zero,
then we
say that it holds {\it quasi everywhere\/} ({\it q.e.\/}) on $E$.
\par
A function $f\colon\om\rightarrow\Rm$ is said to be {\it quasi
continuous\/}
on $\om$ if for every $\e > 0$ there exists a set $E\subseteq\om$ with
${\rm
cap}(E,\om)<\e$ such that the restriction of $f$ to $\om\setminus E$ is
continuous.
A subset $A$ of $\om$ is said to be {\it quasi open\/} if for every $\e >
0$
there exists an open set $A_\e$ with ${\rm cap}(A_\e,\om)<\e$ such
that
$A\cup A_\e$ is an open set.
\par
It is well known (see, for instance, [20]) that for every $u\in\WP$ there
exists a quasi continuous representative of $u$ which is unique up to
sets of
capacity zero, and which is given by
$$
\lim_{r\rightarrow 0^+}\,{1\over{|B_r(x)|}}\int_{B_r(x)}u(y)\,dy
$$
for q.e.\ $x\in\om$.
Throughout this paper we shall use such a quasi continuous
representative to individuate an element of $\WP$. Moreover, we may
also
assume that the quasi continuous representative we are going to choose
is Borel
measurable.
\par
It turns out that for every subset $E$ of $\om$
$$
{\rm cap}(E,\om)\Eq {\rm
inf}\{\|u\|^p_{W^{1,p}(\om)}:u\in W^{1,p}_0(\om)\,, u\Ge 1\
\hbox{q.e.\ on}\ E\}\,.
$$
Actually this infimum is attained by a unique function which is called
the
{\it capacitary potential\/} of $E$. It turns out that this function takes
its
values in $[0,1]$.
\par
A positive Borel measure $\mu$ on $\om$ is said to be {\it absolutely
continuous\/} with respect to capacity if $\mu(B)=0$ whenever
$B\in\B$ and ${\rm cap}(B,\om)=0$.
\medskip
If $E$ is a subset of $\om$ and $F\colon E\rightarrow\Rm$ is a
multivalued
function, i.e., $F$ maps $E$ into the set of all subsets of $\Rm$, then we
say that $F$ is {\it lower semicontinuous at a point $x_0$\/} of $E$ if
for every open subset $G$ of $\Rm$ with $G\cap F(x_0)\Ne\emptyset$
there exists a
neighborhood $U$ of $x_0$ such that for every $y\in U$ we have $G\cap
F(y)\Ne\emptyset$. We say that $F$ is {\it upper semicontinuous in
$x_0$\/} if
for every neighborhood $G$ of $F(x_0)$ there exists a neighborhood $U$
of $x_0$
such that $F(y)\subseteq G$ whenever $y\in U$.
We say that $F$ is {\it quasi lower semicontinuous\/}
(resp. {\it quasi upper semicontinuous\/}) on $\om$ if for every $\e >0$
there
exists a set $E\subseteq\om$ with ${\rm cap}(E,\om)<\e$ such that the
restriction of $F$ to $\om\setminus E$ is lower semicontinuous (resp.
upper
semicontinuous).
 
\medskip
\noindent
{\it Measurability.} Let $(X,\M)$ be a measurable space. If $\mu$ is a
positive
measure on $(X,\M)$ we denote by $\Mu$ the standard
$\mu\!$-completion of $\M$ and we still denote by $\mu$ the
completed measure. If
$\mu$ is $\sigma\!$-finite the $\Mu\!$-measurability is equivalent to
the
$\mu\!$-measurability in the Carath\'eodory sense. Moreover, $\wh
\M$
will denote the universal completion of $\M$, i.e., the intersection
$\cap_\mu\Mu$ for all positive finite measures $\mu$; equivalently,
the
intersection can be extended to all positive  $\sigma\!$-finite measures
$\mu$
(see [10] Ch.III, parag.~4).
\par
It is easy to verify that every quasi continuous function
$u\colon\om\rightarrow\r$ is $\mu\!$-measurable (i.e., ${\cal
B}_\mu\!$-measurable) for every positive Borel measure $\mu$ which
is absolutely
continuous with respect to capacity.
 
\medskip
For convenience, we state
here two results which will play an important role to prove the
measurability of
certain functions. They can be deduced {}from [10], Theorem~III.23 and
Theorem~III.22, respectively.
 
\th{1.1}{{\bf (Projection Theorem)} Let $(X,\M)$ be a measurable space.
If $G$
is an element of  $\M\otimes{\cal B}(\r^d)$, then the projection ${\rm
pr}_X(G)$ belongs to  $\wh \M$.}
 
\th{1.2}{{\bf (Aumann-von Neumann Selection Theorem)} Let $X$ be a
topolo\-gical
space and let $F$ be a multivalued function {}from
$X$ to $\r^d$. If the graph of $F$ belongs to ${\cal
B}(X)\otimes{\cal B}(\r^d)$, then there exists a ${\wh{\cal
B}}(X)$-measurable function which is a selection of $F$ on the set
$\{x\in X: F(x)\Ne\emptyset\}$.}
\medskip
 
\eject
\noindent
{\it Lipschitz projections.} Let us recall some results obtained in
[14], Section 2.
 
\th{1.3}{Let ${\cal C}$ be the family of all non-empty, compact and
convex
subsets of $\Rm$. Then there exists a map $P\colon\Rm\times {\cal
C}\rightarrow\Rm$ satisfying the following properties:
\smallskip
\item{(i)} $P$ is lipschitzian considering on ${\cal C}$ the Hausdorff
metric;
\smallskip
\item{(ii)} $P(\xi,C)\in C$ for every $\xi\in\Rm$, $C\in{\cal C}$,  and
$P(\xi,C)=\xi$ if $\xi\in C$;
\smallskip
\item {(iii)} $d(\xi,C)\Le |\xi - P(\xi,C)|\Le {\sqrt 3}d(\xi,C)$ for every
$\xi\in\Rm$ and $C\in{\cal C}$.
\smallskip}
 
\cor{1.4} {There exists a Lipschitz function
$P_k\colon(\Rm)^{k+1}\rightarrow
\Rm$ satisfying the  following properties:
\smallskip
\item{(i)} $P_k(\xi;\vxi)\in\cox$ for every $\xi$, $\vxi\in\Rm$,
and $P_k(\xi;\vxi)\break =\xi$ if $\xi\in\cox$;
\smallskip
\item {(ii)} $d(\xi,\cox)\Le |\xi - P_k(\xi;\vxi)|\Le {\sqrt 3}d(\xi,\cox)$
for every $\xi$, $\vxi\in\Rm$.
\smallskip}

\rem{1.5}{Let $u$, $\uxi\in\WPL$ and let $C_k(x)=\break\cou$.
If $u(x)\in C_k(x)$ for a.e.\ $x\in\om$, then $u(x)\in C_k(x)$
for q.e.\ $x\in\om$. Indeed, since $P_k$ is lipschitzian, the function
$P_k(u;\uxi)$ is quasi continuous and by assumption
$P_k(u;u_1,\ldots,u_k)= u$ a.e.\ on $\om$. By well-known properties of
Sobolev functions (see [20]) this implies $P_k(u;u_1,\ldots,u_k)= u$
q.e.\ on $\om$, hence $u(x)\in C_k(x)$ for q.e.\ $x\in\om$.}
 
\rem{1.6}{In connection with Theorem~1.3 and Corollary 1.4 we note
that for
every $u\in\WP$, and for every Lipschitz function
$f\colon\Rm\rightarrow\r$
such that $f\circ u\in L^p(\om)$, we have $f\circ u\in W^{1,p}(\om)$
and
$$
|D_i(f\circ u)|\Le L\,|D_iu|\qquad\hbox{a.e.\ on $\om$}\qquad
(i=1,\cdots,n)\,,
$$
where $L$ is the Lipschitz constant of $f$.
This result is classical if $f$ is a $C^1$ function (see, for
instance, [25], Theorem~3.1.9). In the general case it can be obtained by
approximating $f$ with a sequence of $C^1$ functions with Lipschitz
constants bounded by~$L$.}
\eject
 
\parag{2}{A class of ${\bf C^1}\!$-convex functionals: preliminary
properties}
 
Let us first introduce our class $\F$ of $C^1\!$-convex local
functionals.
 
\defin{2.1}{Let $\F$ be the class of all
func\-tio\-nals
$F\colon\WP\times{\cal A}(\om)\rightarrow [0,+\infty]$
satisfying the following properties:
\smallskip
\item{(i)} (lower semicontinuity) for every $A\in{\cal A}(\om)$ the
function
$F(\cdot,A)$ is lower semicontinuous on $\WP$;
\smallskip
\item{(ii)} (measure property) for every $u\in \WP$ the set function
$F(u,\cdot)$ is the trace of a Borel measure on $\om$;
\smallskip
\item{(iii)} (locality property) $F(u,A)= F(v,A)$ whenever $u$,
$v\in\WP$, $A\in {\cal A}(\om)$, and $u|_A = v|_A$;
\smallskip
\item{(iv)} ($C^1\!$-convexity) for every $A\in {\cal A}(\om)$ the
function
$F(\cdot,A)$ is convex on $\WP$ and, in addition, $F(\f u + (1-\f)v,A)\Le
F(u,A) + F(v,A)$ for every $u$, $v\in\WP$ and for every $\f\in
C^1(\om)\cap
W^{1,\infty}(\om)$, with $0\leq\f\leq 1$ on $\om$.
\par}
 
\ex{2.2}{Let $K\colon\om\to\Rm$ be any multifunction with closed
convex values and let
$F\colon\WP\times{\cal A}(\om)\rightarrow [0,+\infty]$ be the
``obstacle functional'' defined by
$$
F(u,A)=\cases{0\,,&if $u(x)\in K(x)$ for q.e.\ $x\in A$,
\cr
\cr+\infty\,,&otherwise.\cr}
$$
Then $F$ satisfies all conditions of Definition~2.1, hence $F\in{\cal F}$.
As mentioned in the introduction, it will be proved in [17] that all
functionals which arise in the study of limits of obstacle problems of the
form (0.1) still belong to the class~$\cal F$.
}
\par
Let $\mu$ and $\nu$ be two positive Borel measures on $\om$ and
let $f\colon\om\times\Rm\to[0,+\infty]$ be a Borel function such that
$f(x,\cdot)$ is convex and lower semicontinuous on $\Rm$ for
$\mu$-a.e.\ $x\in\om$. If $\mu$ is absolutely continuous with respect
to capacity, then the functional
$$
F(u,A)\Eq \int_A f(x, u(x))\,d\mu \, + \, \nu(A)
$$
belongs to the class~$\cal F$. In both examples the lower semicontinuity
follows easily from well known properties of the quasi continuous
representatives of Sobolev functions (see, for instance, [31],
Lemma~2.6.4).

\rem{2.3}{Given a functional $F$ of the class $\cal F$, let us consider the
following extension to $\WP\times\B$:
$$
F(u,B)\Eq \inf\{F(u,A):A\in{\cal A}(\om)\,,B\subseteq A\}\,.\leqno(2.1)
$$
It turns out (see [19], Theorem~5.6) that condition (ii) of Definition~2.1
is equivalent to the assumption that, for every $u\in\WP$, the extension
(2.1) of $F(u,\cdot)$ is a Borel measure on~$\om$.
}
\par
In the sequel, when dealing with Borel sets, we shall always consider
the extension of $F$ given by (2.1).
\par
{}From property (iii) in Definition 2.1 and {}from (2.1) it follows that
$F(u,B)=F(v,B)$ for every $B\in\B$ and for every $u$, $v$ in $\WP$
which coincide in a neighborhood of $B$.
\par
Note that, if $F$ is the obstacle functional of Example~2.2, then, in
general, its extension given by (2.1) does not satisfy
$$
F(u,B)=\cases{0\,,&if $u(x)\in K(x)$ for q.e.\ $x\in B$,
\cr
\cr+\infty\,,&otherwise,\cr}
$$
for every $B\in\B$. For instance, if $n=m=1$ and $\om={]-2,2[}$, let us
consider the obstacle functional
$F\colon\WP\times{\cal A}(\om)\rightarrow [0,+\infty]$ defined by
$$
F(u,A)=\cases{0\,,&if $u(x)\ge x^2$ for q.e.\ $x\in A$,
\cr
\cr+\infty\,,&otherwise.\cr}
$$
Then the extension (2.1) gives $F(1,[0,1])=+\infty$, although
$1\ge x^2$ for every $x\in[0,1]$.
 
\rem{2.4}{If $F$ is a functional of $\F$ and $A\in{\cal A}(\om)$,
then for every finite family $(u_i)_{i\in I}$
of elements of $\WP$ and for every family $(\f^i)_{i\in I}$
of non-negative functions in $\spc$ such that
$\sum_i \f^i=1$ in $\om$, we have
$F(\sum_i\f^i u_i,A)\leq\sum_iF(u_i,A)$.
Indeed, let
$u_1,\ldots,u_r\in\WP$, $\f^1,\ldots,\f^r\in\spc$,
with $\f^i \geq 0$ and $\sum_{i=1}^r\f^i= 1$. It would
be clear, by induction, that
$F(\sum_{i=1}^r\f^iu_i,A)\leq\sum_{i=1}^rF(u_i,A)$
if we had $\f^i >\e$ for every
$i= 1,\ldots,r$ and for a suitable $\e > 0$. Since
$F$ is lower semicontinuous, we can reduce our problem to
this case by considering the coefficients
$\f_\e^i=(\f^i+\e)/(1+r\e)$.
\par
We also notice that, by using the definition (2.1) of $F$ on Borel sets,
property (iv) holds for $A\in\B$, too.}
\bigskip
 
Given $F\in\F$, we now generalize to Borel sets the locality property
(iii) for
$F$ (Proposition~2.6). As a consequence we can single out that part of
the
functional $F$ which is  absolutely continuous with respect to the
capacity
(Proposition~2.8).
\par
For the proof of the locality property on Borel sets we need the
following
remark, which, for future convenience, we state in a slightly more
general form than actually needed here.
 
\lemma{2.5}{Let $s>0$ and $T_s\colon\Rm\rightarrow\Rm$ be the
orthogonal
projection onto the ball ${\overline B}_s(0)$, i.e.,
$$
T_s(\xi)\Eq  {{s}\over{|\xi|\vee{s}}}\,\xi\Eq
\cases{
\xi\,, & if\ \ $|\xi|\leq s\,,$
\cr
{{s}{\xi\over |\xi|}}\,, & if\ \ $|\xi| \geq s\,.$
\cr}
$$
If $u\in\WP$, then $T_s\circ u\in\WPL$ and the sequence
$(T_s\circ u)$ converges, in the strong topology of $\WP$, to $0$  as $s$
goes to $0^+$ and to $u$ as $s$ goes to $+\infty$.}
\bigskip
 
\proof  Since $T_s$ is lipschitzian, for every fixed $u\in\WP$ we have
$T_s\circ u\in\WP$ by Remark 1.6. We prove only the convergence as
$s$ tends to
$0^+$, the other part being analogous. Since $(T_s\circ u)$ converges to
$0$
strongly in $\LP$, there is only to verify the same kind of convergence
for
$(D(T_s\circ u))$. Let $\sigma \geq s$; since orthogonal projections have
Lipschitz constant $1$, the pointwise estimate in Remark 1.6 yields
$$
\displaylines{
\int_{\om}|D(T_s\circ u)|^pdx\ \Eq\
\int_{\{|u|\Le \sigma\}}|D(T_s\circ u)|^pdx +
\int_{\{|u|\G \sigma\}}{|D({s\over\sigma}(T_\sigma\circ u))|^pdx}
\cr
\Le \int_{\{|u|\Le \sigma\}}|Du|^pdx +
({s\over\sigma})^p\,\into{|Du|^pdx}\,.
\cr}
$$
The conclusion now follows taking first the limit as $s$ tends to $0^+$
and
then the limit as $\sigma$ tends to $0^+$.
\endproof
 
\prop{2.6}{{\bf (Locality Property on Borel Sets)}  Let $u$,$v\in\WPL$
and
$B\in\B$.  If $u=v$ q.e.\ on $B$ and $F(u,B)$, $F(v,B) < +\infty$, then
$F(u,B)=F(v,B)$.}
 
\proof {\it Step 1.\ } Assume  $B$ is quasi open. For every $h\in\n$, let
$A_h$
be an open set with ${\rm cap}(A_h,\om)< {1/h}$ and such that $B_h=
B\cup A_h$ is open. Let $w_h$ be the capacitary potential of $A_h$ and
$u_h= u + w_h(v - u)$.
It turns out that $u_h=v$ q.e.\ on $B_h$. Moreover, $(u_h)$ converges to
$u$
in $\WP$ since $(w_h)$ converges to $0$ in
$W^{1,p}(\om)$ and $0\leq w_h\leq 1$ for every $h\in\n$.
\par
Since by assumption $F(u,B) < +\infty$ and $F(v,B) <
+\infty$, it follows that for every given $\e
> 0$ there exist an open set $A$ and a compact set $K$ with $K\subseteq
B\subseteq A$ such that
$F(u,A\setminus K)\L\e$ and
$F(v,A\setminus K)\L\e$.
\par
By the lower semicontinuity of $F$ on open sets we get
$$
F(u,B)\Le F(u,A)\Le\li F(u_h,A)\Le\li \bigl[F(u_h,A\cap B_h) +
F(u_h,A\setminus K)\bigr]\,.\leqno(2.2)
$$
{}From the locality property of $F$ on open sets it follows that
$$
F(u_h,A\cap B_h)\Eq F(v,A\cap B_h)\Le F(v,A)\L F(v,B) +
\e\,.\leqno(2.3)
$$
By approximating $w_h$ in $W^{1,p}(\om)$ with a sequence of
equibounded functions of ${\cal C}^1_0(\om)$, the semicontinuity and
$C^1\!$-convexity of $F$ (properties (i) and (iv)) imply that
$$
F(u_h,A\setminus K)\Le F(u,A\setminus K) + F(v,A\setminus K)\,.
\leqno(2.4)
$$
Hence, $F(u_h,A\setminus K)< 2\e$, and, by (2.2) and (2.3),
$F(u,B)\leq F(v,B) + 3\e$.
Since $\e$ is arbitrary, we can conclude that $F(u,B)\leq F(v,B)$.
Interchanging the roles of $u$ and $v$, we obtain the opposite
inequality.
This proves the theorem when $B$ is quasi open.
 
\medskip\noindent
{\it Step 2.\ } Let now $B$ be a Borel subset of $\om$. For every
$h\in\n$ let
us define $B_h=\break\{x\in \om :|u(x)-v(x)|< 1/h\}$
and
$$
u_h\Eq u + {{1/h}\over{|v-u|\vee (1/h)}}\,(v-u)\Eq u +
T_{1/h}\circ (v-u)\,.
$$
Clearly, $u_h= v$ q.e.\ on $B_h$. By Lemma~2.5 we have the
convergence of
$u_h$ to $u$ in $\WP$. At this point we can introduce the sets $A$ and
$K$ as in Step~1 and proceed in the same way replacing the locality
property
of  $F$
on the open sets with the locality property on the quasi open sets
proved in
Step~1. We have only to remark about the estimate (2.4). Let us notice
that
it is enough to consider $B\subset\subset\om$, hence we
can choose $A\subset\subset\om$; for every
$\om^\p\subset\subset\om$ the
coefficient in the convex combination between $u$ and $v$ defining
$u_h$ can be
approximated in $W^{1,p}(\om^\p)$ by an equibounded sequence of
functions of
$C^1(\overline{\om^\p})$. As $F$ is $C^1\!$-convex and local on
open sets, this suffices to get (2.4) as above.
\endproof
\bigskip
 
Let us consider the function $\nu_0\colon\B\rightarrow [0,+\infty]$
defined by
$$
\nu_0(B)\Eq {\Inf\{F(v,B): v\in\WPL\}}\leqno(2.5)
$$
for every $B\in\B$. Moreover, for every $A\in\A$ we  define
$$
\d F(\cdot,A)=\{u\in\WP:F(u,A)< +\infty\}\,.
$$
Then the following proposition holds.
 
\prop {2.7} {For every $\om'\in\A$ with $\d
F(\cdot,\om')\cap\LL\Ne\emptyset$, the restriction of $\nu_0$ to
${\cal B}(\om')$
is a positive finite Borel measure.}
 
\proof It is clear that $\nu_0$ is an increasing function,
$\nu_0(\emptyset)=0$ and, in view of the definition of $F(u,\cdot)$ on
$\B$,
that $\nu_0(B)=\inf\{\nu_0(A):A\in{\cal A}(\om'), A\supseteq B\}$.
Hence, by
Proposition~5.5 and Theorem~5.6 in [19], we have only to prove that
$\nu_0$ is superadditive, subadditive, and inner regular on
${\cal A}(\om')$. The superadditivity
comes immediately {}from the definition of $\nu_0$ and {}from the
additivity of
$F$ in the second variable. Let us now prove that for every $A_1$,
$A_2$,
$A_2^\p\in{\cal A}(\om')$ with $A_2^\p\subset\subset A_2$ we have
$$
\nu_0(A_1\cup A_2^\p)\Le \nu_0(A_1) + \nu_0(A_2)\,.\leqno(2.6)
$$
We can assume that $\nu_0(A_1) + \nu_0(A_2)< +\infty$. Then,
for every $\e > 0$ there exist two functions
$u_1$, $u_2$ in $\WPL$ such that
$$
\nu_0(A_1) + {\e\over 2}\G
F(u_1,A_1)\qquad\qquad \nu_0(A_2)  + {\e\over 2}\G F(u_2,A_2)\,.
$$
Let $\f\in C^1_0(A_2)$, with $\f= 1$ on a neighborhood of
${\overline {A_2^\p}}$ and $0\leq \f \leq 1$.
We set $u=(1-\f)u_1 +\f u_2$. By Remark 2.3 it follows that
$$
\centerline{\vbox{\halign{\hfil$\displaystyle#
$&$\displaystyle#$\hfil\cr
\nu_0(A_1\cup A_2^\p)&\Le F(u,A_1\cup A_2^\p)
\cr
&\Le F(u_1,A_1\setminus A_2) + F(u_2,{\overline {A_2^\p}}) +
F((1-\f)u_1+\f u_2,(A_2\setminus{\overline{A_2^\p}})\cap A_1)\,.
\cr}}}
$$
Furthermore, the $C^1\!$-convexity of $F$ permits to estimate the last
term in
the above inequality by $F(u_1,(A_2\setminus{\overline{A_2^\p}})\cap
A_1) +
F(u_2,(A_2\setminus{\overline{A_2^\p}})\cap A_1)$. Hence,
$$
\nu_0(A_1\cup A_2^\p)\Le F(u_1,A_1) + F(u_2,A_2) \L
\nu_0(A_1) + \nu_0(A_2) + \e\,.
$$
Thus we obtain (2.6). This inequality will give the subadditivity of
$\nu_0$
once inner regularity will be proved.
\par
Since $\d F(\cdot,\om')\cap\LL\Ne\emptyset$, we can find $u\in\WPL$
such that
$F(u,\om')<+\infty$. Therefore, given $A\in{\cal A}(\om')$ and $\e>0$
there exists
$A^{\pp}\in{\cal A}(\om')$ with $A^{\pp}\subset\subset A$ and
$F(u,A\setminus{\overline {A^{\pp}}})\leq\e$; it follows that
$\nu_0(A\setminus{\overline {A^{\pp}}})\leq\e$. Let $A^\p\in{\cal
A}(\om')$ such
that $A^{\pp}\subset\subset A^\p\subset\subset A$. By (2.6) we have
$$
\nu_0(A)\Le \nu_0(A^\p) +J\nu_0(A\setminus{\overline
{A^{\pp}}})\leq
\nu_0(A^\p) +\e\,. $$
We conclude that $\nu_0(A)=\sup\{\nu_0(A^\p):A^\p\in{\cal A}(\om'),
A^\p\subset\subset A\}$, i.e., the inner regularity of $\nu_0$.
\endproof
 
\prop {2.8} {For every $A\in {\cal A}(\om)$ and $u\in\d
F(\cdot,A)\cap\LL$,
the function
$F(u,\cdot)-\nu_0(\cdot)$ is a positive Borel measure on $A$
which is absolutely continuous with respect to capacity.}
 
\proof Let $u\in\WPL$ such
that $F(u,A)< +\infty$. As $\nu_0$ is a finite Borel measure on $A$
(Proposition~2.7) and $\nu_0(\cdot)\Le F(u,\cdot)$ by (2.5), we
conclude that
$F(u,\cdot)-\nu_0(\cdot)$ is a positive Borel measure on $A$. Let us fix
$B\in{\cal B}(A)$ with ${\rm cap}(B,\om)= 0$. For every $v\in\WPL$
with
$F(v,B)< +\infty$, we have $v= u$ q.e.\ on $B$; hence, by Proposition~2.6
we
conclude that $F(v,B)=F(u,B)$. Since
$$ \nu_0(B)\Eq {\Inf}\{F(v,B):
v\in\WPL,\ F(v,B)\L+\infty\}
$$
it follows that
$\nu_0(B)= F(u,B)$, i.e., $F(u,B)-\nu_0(B)= 0$.
\endproof
\bigskip
 
We now conclude this section by giving a basic estimate for $F$ on the
convex
hull of a finite  number of functions in $\WPL$.
 
\prop{2.9}{Let $u$,$\uxi\in\WPL$. Assume that $u(x)\in\cou$
for a.e.\ $x\in\om$. Then
$$
F(u,B)\Le\s F(u_i,B)
$$
for every Borel set $B$ in $\om$.}
 
\proof  In view of the definition of $F$ on Borel sets, it is enough to
prove
the inequality for every open set $B$ with $B\subset\subset\om$.
Hence, let
$B$ be such a set. By means of the Density Theorem~2.9 in~[14], we can
easily find a sequence of functions
$\f_h\colon\r^n\rightarrow\Sigma_k$ such
that $\f_h\in C^\infty(\r^n,\r^k)$ and
$$\s\f^i_h u_i\ {\mapright h}\ \ u$$
strongly in $W^{1,p}(A,\Rm)$, where $A$ is a neighborhood of
$\overline B$.
Then, {}from the lower semicontinuity of  $F$ and the locality property
on open
sets we obtain
$$
F(u,B)\Le\li F(\s\f^i_h u_i,B)\,.
$$
The conclusion follows now {}from the $C^1\!$-convexity of $F$ and
{}from
Remark 2.4.
\endproof
\bigskip
 
\parag{3}{Integral representation on moving polytopes}
 
The aim of this section is the integral representation of the functionals
in
$\F$ when restricted to the pointwise convex hull of a finite number of
functions
in $\WPL$ (see Theorem~3.7).
\par
Let $F\in\F$, $k\in\n$, and $u_1,\ldots,u_k\in\WPL$ be fixed.
Throughout
this section we assume that $F(u_i,\om)< +\infty$ for every
$i=1,\ldots,k$. We
point out that our proof first produces a kind of integral representation
of
$F$ on the functions of the form $u=\s\psi^i u_i$, where
$\psi\colon\om\to\Sigma_k$, with the integrand depending on the
coefficient
$\psi$ rather than on the function $u$ itself (Theorem~3.3). For a
constant
$\psi$, this result is contained in the following lemma.
 
\bigskip
Let $\nu_0$ be the set function introduced in (2.5). Under the present
assumptions, Proposition~2.7 tells us that $\nu_0$ is a finite Borel
measure on
$\om$. Let $\mu$ be a positive finite Borel measure on $\om$ with
${\rm supp}\,\mu=\om$ such
that $\mu$ is absolutely continuous with respect to
capacity and $\mu (\cdot)\geq\s\bigl(F(u_i,\cdot)-\nu_0(\cdot)\bigr)$
(in
view of Proposition~2.8, such a  $\mu$ can be obtained, for example, by
adding to
$\s\bigl(F(u_i,\cdot)-\nu_0(\cdot)\bigr)$ the positive measure $fd{\cal
L}^n$,
where $f\in L^1(\om)$, $f > 0$ on $\om$, and ${\cal L}^n$ is the
$n$-dimensional Lebesgue measure).

\lemma {3.1}{ For every $x\in\om$ and $\lambda\in\Sigma_k$ we
define
$u_\lambda (x) = \s\lambda^iu_i(x)$ and
$$
g(x,\lambda)\Eq\lsr{{F(u_\lambda,B_r(x)) -J\nu_0(B_r(x))}\over
\mu(B_r(x))}\,.\leqno(3.1)
$$
Then
\item {(i)} for every $x\in\om$ the function $g(x,\cdot)$
is convex and continuous in $\Sigma_k$;
\smallskip
\item {(ii)} for every $\lambda\in\Sigma_k$ the
function $g(\cdot,\lambda)$ is Borel measurable on $\om$;
\smallskip
\item {(iii)}
$F(u_\lambda,B)=\intb{g(x,\lambda)\,d\mu} + \nu_0(B)$ for every
$\lambda\in\ris$ and $B\in\B$.
\smallskip}
 
\proof  {}From the definition of $\mu$ and $\nu_0$ and {}from the
convexity of
$F$, it follows immediately that $0\leq g\leq 1$ on $\om\times\ris$,
and that
$g(x,\cdot)$ is convex on $\ris$ for every $x\in\om$. Then, by
Theorem~10.3
in~[28], for every $x\in\om$ the function $g(x,\cdot)$ can be extended
in one and
only one way to a continuous convex function, still denoted by $g$, on
the whole
of $\Sigma_k$. Hence, $0\leq g\leq 1$ on $\om\times\Sigma_k$ and (i)
holds. \par
Let us proof (ii). If $\alpha$ is a positive Borel measure on $\om$,
the function $r\mapsto \alpha(B_r(x))$ is left continuous for every
$x\in\om$. This implies that the upper limit which appears in (3.1) can
equivalently be taken as $r\rightarrow 0^+$ with $r\in\q$. Moreover,
the
function $x\mapsto \alpha(B_r(x))$ is lower semicontinuous for every
$r>0$ and hence Borel measurable, too. It follows that the function
$g(\cdot,\lambda)$ is
Borel measurable for every $\lambda\in\ris$. For
$\lambda\in\partial\Sigma_k$,
$g(\cdot,\lambda)$ is the pointwise limit of a sequence
$g(\cdot,\lambda_n)$
with $\lambda_n\in\ris$; therefore, $g(\cdot,\lambda)$ is Borel
measurable on
$\om$ for every $\lambda\in\Sigma_k$.
\par
By the Besicovitch differentiation theorem (see, e.g., [31], Section 1.3),
we
have
$$ F(u_\lambda,B)\Eq\intb{g(x,\lambda)\,d\mu} +
\nu_0(B)
$$
for every $B\in\B$ and for every $\lambda\in\ris$.
\endproof
\bigskip
 
Before extending the previous result to non constant $\lambda$'s, we
observe that the following selection lemma holds.
 
\lemma{3.2}{Let $u\in\WPL$ such that $u(x)\in\cou$ for a.e.\
$x\in\om$.
Then, there exists a $\wh{\cal B}(\om)$-measurable function
$\psi\colon\om\to\Sigma_k$
such that
$$
u(x)\Eq\s\psi^i(x)u_i(x)\qquad\hbox{for q.e.\ $x\in\om$.}\leqno(3.2)
$$}
\vskip-1.2truecm\noindent
\proof Let us fix some
quasi continuous Borel measurable representatives of $u,u_1,\ldots,u_k$
(see
Section 1). Let $\Lambda(x)=\bigl\{\lambda\in\Sigma_k :
u(x)=\s\lambda^i u_i(x)\bigr\}$ for every $x\in\om$. $\Lambda$
is a multivalued function {}from $\om$ to $\Sigma_k$ with non-empty
closed values
for q.e.\ $x\in\om$ (see Remark~1.5). It is clear that
graph$\Lambda\in\B\otimes{\cal B}(\r^k)$.
\par
By the Aumann-von Neumann Selection Theorem~1.2 there exists a
${\wh{\cal B}}(\om)$-measur\-able selection~$\psi$
of the multivalued function $\Lambda$ and, by the definition of
$\Lambda (x)$,
the function $\psi$ satisfies (3.2).
\endproof
 
\th{3.3}{ Let $u\in\WPL$ with $u(x)\in\cou$ for a.e.\ $x\in\om$. If
$g$
is the function given by Lemma~3.1 and
$\psi\colon\om\rightarrow\Sigma_k$
is a ${\wh{\cal B}}(\om)$-measurable function such that (3.2) holds,
then $g(\cdot,\psi(\cdot))$ is $\mu\!$-measurable and
$$
F(u,A)\Eq\int_A{g(x,\psi(x))\,d\mu} + \nu_0(A)\leqno (3.3)
$$
for every $A\in\A$.}
\medskip
 
Let us explicitly notice that if $u$ is as above, then $F(u,\om)<+\infty$
by
Proposition~2.9.
The proof of Theorem~3.3 heavily relies on the following
approximation lemma, which essentially reduces the problem to the case
of a
constant $\psi$.
 
\lemma{3.4}{Let $u$ and $\psi$ be as in Theorem~3.3. Let
$\lambda\in\ris$ with
$d(\lambda,\partial\Sigma_k) = \eta > 0$, let $0<\e <\eta$ and $B\in\B$
such
that $|\psi(x) - \lambda|\leq\e$ for q.e.\ $x\in B$.
Then,
$$
|F(u,B) - F(u_\lambda,B)|\Le {\e\over\eta}\s F(u_i,B)\,.
$$
}
 
\proof Let us define
$v= u + t(u_\lambda - u)$ on $\om$,
with $t= 1 + \eta/\e$. It turns out that $v(x)\in
\cou$ for q.e.\ $x\in B$ and
$$
u_\lambda\Eq {1\over t}v + (1-{1\over t})u\qquad {\rm on}\ \om\,.
$$
In order to get {}from $v$ a function which belongs to ${\rm
co}\{u_1,\ldots,u_k\}$ a.e.\ on $\om$, we consider the projection
$w= P_k(v;\uxi)$ as defined in Corollary 1.4. Set
$$
z\Eq {1\over t}w + (1-{1\over t})u\qquad {\rm on}\ \om\,.
$$
By Remark~1.6, the function $w$, and hence $z$, belongs to $\WPL$.
Moreover,
since $w= P_k(v;u_1,\ldots,u_k)= v$ q.e.\ on $B$ (see Corollary
1.4), we have $z= u _\lambda$ q.e.\ on $B$. By the convexity of $F$
and Proposition~2.9
$$
\centerline{\vbox{\halign{\hfil$\displaystyle#
$&$\displaystyle#$\hfil\cr F(z,B)&\Le {1\over t}F(w,B) + (1-{1\over
t})F(u,B)
\cr
&\Le {\e\over\eta}\s F(u_i,B) + F(u,B)\,.
\cr}}}
$$
In view of the locality property of $F$ on Borel sets (Proposition~2.6) we
have
$F(z,B)= F(u_\lambda,B)$; thus
$$
F(u_\lambda,B)\Le F(u,B) + {\e\over\eta}\s F(u_i,B)\,.
$$
The inequality
$$
F(u,B)\Le F(u_\lambda,B) + {\e\over\eta}\s F(u_i,B)
$$
can be obtained analogously defining now $v= u_\lambda + t(u-
u_\lambda)$
with $t= \eta/\e > 1$.
\endproof
\bigskip
 
\noindent
{\it Proof of Theorem~3.3.}\ Let us fix
$A\in{\cal A}(\om)$.
\smallskip
\noindent
{\it Step 1.} Assume that $\psi(x)\in\ris$ for every $x\in\om$.
\par\noindent
Given $\eta > 0$, let us define
$\Sigma_{k,\eta}= \bigl\{\lambda\in\Sigma_k:
d(\lambda,\partial\Sigma_k) \geq \eta\bigr\}$
and ${A_\eta =\psi^{-1}\bigl(\Sigma_{k,\eta}\bigr)\cap A}$.
For every  $\e\in ]0,\eta[$ we can fix a finite partition $(B_j)_{j\in J}$
of $\Sigma_{k,\eta}$ by means of Borel sets having diameter less than
$\e$, and
a family $(\lambda_j)_{j\in J}$ of elements of $\Sigma_{k,\eta}$ such
that
$\lambda_j\in B_j$ for every $j\in J$. Let us define $E_j
=\psi^{-1}\bigl(B_j)\cap A$ for every $j\in J$.  Since $\psi$
is ${\wh {\cal B}}(\om)$-measurable, the sets $A_\eta$ and $E_j$ are in
${\wh {\cal B}}(\om)$. According to the convention made in Section 1,
for every $z\in\WP$ the completion of the measures
$\mu$, $\nu_0$ and $F(z,\cdot)$ will be still denoted by $\mu$,
$\nu_0$ and  $F(z,\cdot)$ . By Lemma~3.4, for every $j\in J$
$$
|F(u,E_j)-F(u_{\lambda_j},E_j)|\Le {\e\over\eta}\sum_{i=1}^{k}
F(u_i,E_j)\,.
$$
This and Lemma~3.1 imply
$$
\centerline{\vbox{\halign{\hfil$\displaystyle#
$&$\displaystyle#$\hfil\cr
F(u,A_\eta)-\nu_0(A_\eta)&\Eq
\sum_{j\in J}\bigl[F(u,E_j)-\nu_0(E_j)\bigr]
\cr
&\Le \sum_{j\in J}\bigl[F(u_{\lambda_j},E_j)-\nu_0(E_j) +
{\e\over\eta}\s F(u_i,E_j)\bigr]
\cr
&\Eq\sum_{j\in J}\bigl[\int_{E_j}g(x,\lambda_j)\,d\mu
+{\e\over\eta}\s F(u_i,E_j)\bigr]\,.
\cr}}}
$$
Since $g(x,\cdot)$ is convex and bounded by 1 in $\Sigma_k$, it is
Lipschitz
continuous in $\Sigma_{k,\eta}$ with constant $1/\eta$; thus
$$
\eqalign{F(u,A_\eta)-\nu_0(A_\eta)&\Le
\sum_{j\in J}\bigl[\int_{E_j}g(x,\psi(x))\,d\mu
+{1\over\eta}\int_{E_j}|\psi(x)-\lambda_j|\,d\mu
+{\e\over\eta}\s F(u_i,E_j)\bigr]
\cr
&\Le\int_{A_\eta}g(x,\psi(x))\,d\mu +
{\e\over\eta}\bigr[\mu(A_\eta)+\s F(u_i,A_\eta)\bigr]\,.
\cr}
$$
(Note that $g(\cdot,\psi(\cdot))$ is $\mu$-measurable since $g$ is
Borel measurable and $\psi$ is
${\wh{\cal B}}(\om)$-measurable).
Since $\e$ is arbitrary, we get
$$
F(u,A_\eta)-\nu_0(A_\eta)\Le\int_{A_\eta}g(x,\psi(x))\,d\mu\,.
$$
Now, taking into account that $\psi(x)\in\ris$ for every $x\in\om$ and
letting $\eta\rightarrow 0^+$, we obtain
$$
F(u,A) - \nu_0(A)\Le\int_A g(x,\psi(x))\,d\mu\,.
$$
\smallskip
\noindent
The reverse inequality can be obtained in a completely analogous way.
 
\bigskip
\noindent
{\it Step 2.}\  We consider now the general case
$\psi\colon\om\rightarrow\Sigma_k$.
\par\noindent
Let $b_0= {1\over k}(e_1 +\cdots+e_k)$ be the
barycenter of $\Sigma_k$ ($e_1,\ldots,e_k$ are the elements of the
standard
basis of $\r^k$). For every $0\leq \sigma \leq 1$ define $\psi_\sigma=
b_0 +
\sigma(\psi - b_0)$ and
$u_\sigma=\s\psi_\sigma^iu_i= u_0 + \sigma(u-u_0)$, where
$u_0={1\over
k}(u_1+\cdots+u_k)$.
If $0\leq\sigma<1$ then $\psi_\sigma(x)\in\ris$ for every
$x\in\om$ (see~[28], Theorem~6.1); therefore, by Step 1 we have
$$
F(u_\sigma,A)\Eq\int_A g(x,\psi_\sigma(x))\,d\mu +\nu_0(A)\,.
$$
Observe now that the lower semicontinuity and the convexity of $F$
imply that
$$
\lim_{\sigma\rightarrow 1^-}F(u_\sigma,A)\Eq F(u,A)\,.
$$
Moreover, the continuity of $g(x,\cdot)$ and the dominated convergence
theorem
yield
$$
\lim_{\sigma\rightarrow 1^-}\int_A g(x,\psi_\sigma(x))\,d\mu\Eq
\int_A g(x,\psi(x))\,d\mu\,.
$$
This concludes the proof.
\endproof
\bigskip
 
We point out that the values $u(x)$ of the function $u$ enter the
integral
representation of Theorem~3.3 through the parameters $\psi^i(x)$ for
which
$u(x)=\sum_{i=1}^k\psi^i(x)u_i(x)$. When looking for an integrand
depending directly on the values of $u$, the main difficulty we meet is
that
the expression of $u(x)$ as a convex combination of
$u_1(x),\ldots,u_k(x)$ is
not necessarily unique. This problem is essentially overcome by the
following
lemma.
 
\lemma {3.5}{ Let $u\in\WPL$ and $u(x)\in\cou$ for a.e.\ $x\in\om$.
Let $g$ be
the function defined in Lemma~3.1 and
$$
N\Eq \bigl\{x\in\om :\exists\lambda,\lambda^\p\in\Sigma_k\ \,
u(x)\Eq u_\lambda(x)\Eq u_{\lambda^\p}(x)
\  and\  g(x,\lambda)\Ne g(x,\lambda^\p)\bigr\}
$$
(N is defined up to sets of zero capacity). Then, $N\in{\wh{\cal
B}}(\om)$ and
$\mu(N)= 0$.}
 
\proof  As in the proof of Lemma~3.2 we can fix quasi continuous, Borel
measurable representatives of $u$, $u_1,\ldots, u_k$; the set $N$ is now
well
defined all over $\om$.  Consider the multivalued map $\Gamma$
{}from $\om$ into $\Sigma_k\times\Sigma_k$ defined by
$$
\Gamma
(x)\Eq\bigl\{(\lambda,\lambda^\p)\in\Sigma_k\times\Sigma_k :
u(x)\Eq u_\lambda(x)\Eq u_{\lambda^\p}(x)
\ {\rm and}\ g(x,\lambda)\Ne g(x,\lambda^\p)\bigr\}\,.
$$
Arguing as in the proof of Lemma~3.2 and taking into account that $g$ is
Borel measurable on $\om\times\Sigma_k$, we
obtain that ${\rm graph}\Gamma\in\B\otimes{\cal B}(\r^k)\otimes
{\cal B}(\r^k)$.
By the Projection Theorem (Theorem~1.1) we
get $N\in{\wh {\cal B}}(\om)$. The Aumann-von Neumann Selection
Theorem
(Theorem~1.2) implies the existence of two ${\wh {\cal
B}}(\om)$-measurable functions $\sigma_1$, $\sigma_2\colon
\om\rightarrow\Sigma_k$ such that $(\sigma_1|_N,\sigma_2|_N)$ is a
selection of
$\Gamma$ on $N$.
\noindent
Define for $j= 1,2$
$$
\psi_j\Eq
\cases{
\psi\,, & on $\om\setminus N\,,$
\cr
\sigma_j\,, & ${\rm on}\  N\,,$
\cr}
$$
where $\psi\colon\om\rightarrow\Sigma_k$ is the
${\wh {\cal B}}(\om)$-measurable function given in Lemma~3.2. Then
$\psi_1$ and
$\psi_2$ are ${\wh {\cal B}}(\om)$-measurable functions such that
$$
u(x)\Eq\sum_{i=1}^k \psi_1^i(x)u_i(x)\Eq\sum_{i=1}^k \psi_2^i(x)u_i(x)\
\ \
{\rm for\ \ q.e.\ } x\in\om\,,\leqno (3.4)
$$
$$
g(x,\psi_1(x))\Ne g(x,\psi_2(x))\ \ \ {\rm for}\ {\rm every}\  x\in
N\,.\leqno
(3.5)
$$
By Theorem~3.3, (3.4) implies that
$$
\int_A {g(x,\psi_1(x))\,d\mu}\Eq \int_A {g(x,\psi_2(x))\,d\mu}
$$
for every $A\in\A$. Hence,
$$
g(x,\psi_1(x))\Eq g(x,\psi_2(x))\qquad\hbox{for $\mu$-a.e.\
$x\in\om\,.$}
$$
Together with (3.5) this yields that $\mu (N)\Eq 0$.
\endproof
\bigskip
 
For future convenience we single out a technical remark about
measurability of
functions.
 
\rem{3.6}{(i)\ \ Let $T$ be a Borel subset of $\Rm$. Given
$g\colon\om\times\Rm\times T\rightarrow [0,+\infty]$, let us define
$f\colon\om\times\Rm\rightarrow [0,+\infty]$ by setting
$f(x,\xi)=\inf_{t\in
T}g(x,\xi,t)$. If $g(x,\cdot,t)$ is continuous for every $(x,t)\in\om\times
T$,
uniformly with respect to $t\in T$, then $f(x,\cdot)$ is continuous for
every
$x\in\om$. Assume, in addition, that $g(\cdot,\xi,\cdot)$ is ${\wh{\cal
B}}(\om)\otimes{\cal B}(\Rm)$-measurable for every $\xi\in\Rm$. Then
$f(\cdot,\xi)$ is ${\wh{\cal B}}(\om)$-measurable for every
$\xi\in\Rm$, hence
$f$ is  ${\wh{\cal B}}(\om)\otimes{\cal B}(\Rm)$-measurable.
\par
Indeed, given $\xi\in\Rm$ and $s\in\r$, the set
$$E_s=\{x\in\om:f(x,\xi)<s\}=\{x\in\om:\exists t\in T\ \,g(x,\xi,t)<s\}$$
is the projection on $\om$ of the set $\{(x,t)\in\om\times T:
g(x,\xi,t)<s\}$.
Since $(x,t)\mapsto g(x,\xi,t)$ is
${\wh{\cal B}}(\om)\otimes{\cal B}(\Rm)$-measurable, by the
Projection
Theorem (Theorem~1.1) we get $E_s\in {\wh{\cal B}}(\om)$.
\medskip
\noindent
(ii)\ \ Let $f\colon\om\times\Rm\rightarrow [0,+\infty]$ be a
${\wh{\cal B}}(\om)\otimes{\cal B}(\Rm)$-measurable function. Then
for every
positive finite Borel measure $\mu$ on $\om$ there exists a set
$N\in\B$ with
$\mu(N)=0$ such that $f|_{(\om\setminus N)\times\Rm}$ is a Borel
function.
\par
Indeed, for every $E\in {\cal B}_\mu(\om) \otimes{\cal B}(\Rm)$ there
exists
$N\in\B$ with $\mu(N)=0$ such that
$E\setminus(N\times\Rm)\in\B\otimes{\cal
B}(\Rm)$.
}
\bigskip
 
Finally, we are able to prove the main result of this section, i.e., the
integral representation on the pointwise convex combinations of a finite
number
of fixed functions.
 
\th{3.7}{Let $F\in\F$, $k\in\n$, $u_1,\ldots,u_k\in\WPL$, and assume
$F(u_i,\om)< +\infty$ for every $i= 1,\ldots,k$. Let $\nu_0$ be the
positive
finite Borel measure introduced in Proposition~2.7. Then, there exist a
positive
finite Borel measure $\mu$ on $\om$, absolutely continuous with
respect to
capacity, and a function $f\colon\om\times\Rm\rightarrow [0,+\infty]$
with the
following properties:
\smallskip
\item{(i)} for every $x\in\om$ the function $f(x,\cdot)$
is convex and lower semicontinuous on $\Rm$;
\smallskip
\item{(ii)} f is ${\wh{\cal B}}(\om)\otimes{\cal B}(\Rm)$-measurable;
\smallskip
\item {(iii)} for every $u\in\WPL$ such that $u(x)\in\cou$ for
a.e.\ $x\in\om$, the function $f(\cdot,u(\cdot))$ is
$\mu$-measurable on $\om$ and $F(u,A)=\int_A {f(x,u(x))\,d\mu} +
\nu_0(A)$ for
every  $A\in\A$.
Moreover, the restriction of $f(x,\cdot)$ to $\cou$ is
continuous for $\mu$-a.e.\ {$x\in\om$}.
\smallskip
}
 
\proof  Let us fix quasi continuous Borel measurable representatives of
$u_1,\ldots,u_k$. For every $\lambda\in\Sigma_k$ and $x\in\om$ we
set
$u_\lambda (x) =
\sum_{i=1}^k \lambda^i u_i(x)$ and $C_k(x) =
{\rm co}\{u_1(x),\ldots,u_k(x)\}$. Let us define the
function $f\colon\om\times\Rm\to [0,+\infty]$ as
$$
f(x,\xi) =
\cases{
\displaystyle\Inf_{\scriptstyle\lambda\in\Sigma_k\atop\scriptstyle{
u_\lambda (x) =\xi}}g(x,\lambda)\,, & if $\ \xi\in C_k(x)\,,$
\cr
\Extra
\cr
+\infty\,, & otherwise$\,,$
\cr}
$$
where $g$ is the function introduced in Lemma~3.1.
\par
Let us prove (i). Fix $x\in\om$; the multivalued function {}from
$C_k(x)$ to
$\Sigma_k$ defined by $\xi\mapsto\{\lambda\in\Sigma_k:
u_\lambda(x)=\xi\}$ has
closed graph and compact range; hence it is
upper semicontinuous.
Moreover, by the continuity of $g(x,\cdot)$, the function
${S\mapsto{\Inf}\{g(x,\lambda)\colon\lambda\in S\}}$, defined on the
compact
subsets of $\Sigma_k$, is continuous with respect to the Hausdorff
metric and
decreasing with respect to inclusion. Therefore, we can deduce the lower
semicontinuity of $f(x,\cdot)$ on $C_k(x)$. This immediately implies the
lower
semicontinuity of $f(x,\cdot)$ on $\Rm$, while the convexity of
$f(x,\cdot)$ can
be easily verified directly. Hence, (i) holds true.
\par
Let us prove (ii). For every $x\in\om$ let us consider the Moreau-Yosida
transforms of $f(x,\cdot)$, defined by
$$
f_s(x,\xi)\Eq\Inf_{\eta\in C_k(x)}[f(x,\eta) + s|\xi - \eta|]
\qquad (s\in\n)
$$
for every $\xi\in\Rm$. Since $f(x,\cdot)$ is lower semicontinuous on
$\Rm$,
for every $x\in\om$ and $\xi\in\Rm$ we have
$$
f(x,\xi)= \Sup_{s\in\n}f_s(x,\xi)\,.\leqno(3.6)
$$
Let us prove that for every $\xi\in\Rm$, $f_s(\cdot,\xi)$ is ${\wh{\cal
B}}(\om)$-measurable. Note that
$$
\centerline{\vbox{\halign{\hfil$\displaystyle#
$&$\displaystyle#$\hfil\cr
f_s(x,\xi)&=\Inf_{\eta\in C_k(x)}\Inf_{\lambda\in\Sigma_k\atop
u_\lambda (x)=\eta} [g(x,\lambda ) + s|\xi - \eta|]
\cr
&= \Inf_{\lambda\in\Sigma_k}[g(x,\lambda) + s|\xi -
u_\lambda (x)|]\,.\cr}}}
$$
Remark~3.6 shows that $f_s$ is ${\wh{\cal B}}(\om)\otimes{\cal
B}(\Rm)$-measurable; by (3.6) the same is true for $f$.
\par
Let us now turn to the proof of (iii). Fix $u\in\WPL$ such that
$u(x)\in\cou$
for a.e.\ $x\in\om$ and choose a quasi continuous representative
of $u$. Since such a representative is ${\cal
B}_\mu(\om)$-measurable (recall that $\mu$ is absolutely continuous
with respect
to capacity), $f(\cdot,u(\cdot))$ is $\mu$-measurable on $\om$. By
Lemma~3.2
and Theorem~3.3, for every $A\in\A$ we have
$$
F(u,A)\Eq\int_A {g(x,\psi(x))\,d\mu} + \nu_0(A)\,,
$$
where  $\psi\colon\om\rightarrow\Sigma_k$ is a
${\wh {\cal B}}(\om)$-measurable function such that (3.2) holds. By the
definition of $f$ we have
$$
f(x,u(x))\Eq \Inf_{\lambda\in\Sigma_k\atop {u_\lambda(x)\Eq  u(x)}}
g(x,\lambda)\qquad\hbox{for q.e.\ $x\in\om\,.$}
$$
Let $N$ be the set given in Lemma~3.5. Then
$$
f(x,u(x))\Eq g(x,\psi(x))\qquad {\rm for\ q.e.\ }\
x\in\om\setminus N\,.
$$
This proves (iii) since $\mu(N)=0$.
\endproof
\bigskip
 
The following proposition shows that, given the measures $\mu$ and
$\nu$, the
function $f$ obtained in the integral representation theorem is
essentially
unique.
 
\prop{3.8} {Let $F$, $u_1,\ldots,u_k$ be as in Theorem~3.7. Let $\mu$
and
$\nu$ be two positive finite Borel measures on $\om$, with $\mu$
absolutely
continuous with respect to capacity. Assume that two functions $f_1$,
$f_2\colon\om\times\Rm\rightarrow [0,+\infty]$ satisfy conditions
(i)--(iii) of Theorem~3.7 with $\nu_0$ replaced by $\nu$. Then
$f_1(x,\xi)=f_2(x,\xi)$ for $\mu$-a.e.\ $x\in\om$ and for every
$\xi\in\cou$.}
 
\proof
{}From the finiteness of $F(u_i,\om)$ for $i=1,\ldots,k$ and {}from
property
(iii), we deduce that $f_1(\cdot,u_i(\cdot))<+\infty$,
$f_2(\cdot,u_i(\cdot))<+\infty$ $\mu$-a.e.\ on $\om$. The convexity of
$f_1(x,\cdot)$ and $f_2(x,\cdot)$ then guarantees that $f_1(x,\cdot)$
and
$f_2(x,\cdot)$ are finite on $\cou$ for $\mu$-a.e.\ $x\in\om$.
By Theorem~10.1 in~[28], it follows that $f_1(x,\cdot)$ and
$f_2(x,\cdot)$ restricted to ${\rm ri\,}\cou$ are continuous for
$\mu$-a.e.\ $x\in\om$.
By (iii), for every
$\lambda\in\Sigma_k\cap\q^k$ we have
$$
\int_A f_1(x,u_\lambda(x))\,d\mu = \int_A
f_2(x,u_\lambda(x))\,d\mu
$$
for every $A\in\A$; hence, there exists a set $N\in\B$,
$\mu(N)= 0$ such that $f_1(x,u_\lambda(x))= f_2(x,u_\lambda(x))$
for every $x\in\om\setminus N$ and $\lambda\in\Sigma_k\cap\q^k$.
Since
the functions
$f_1(x,\cdot)$ and $f_2(x,\cdot)$ restricted to ${\rm
ri~co}\{u_1(x),\ldots,u_k(x)\}$ are continuous,
we have $f_1(x,\xi)= f_2(x,\xi)$
for
every $x\in\om\setminus N$ and $\xi\in{\rm
ri~co}\{u_1(x),\ldots,u_k(x)\}$.
By the continuity along line segments
([28], Corollary 7.5.1) $f_1(x,\xi)= f_2(x,\xi)$ for every $\xi\in{\rm
co}\{u_1(x),\ldots,u_k(x)\}$.
\endproof
\bigskip
 
\parag{4}{Auxiliary lemmas}
 
We collect here some results we shall use in the next section.
 
\lemma{4.1}{Let $X$ be a separable metric space and let $F\colon
X\rightarrow\r$ be lower semicontinuous. Then there exists a countable
subset $D$ of $X$ with the following property: for every $x\in X$ there
exists a sequence $(x_h)$ in $D$ converging to $x$ and such that
$(F(x_h))$
converges to~$F(x)$.}
 
\proof  It is enough to take a countable dense subset $E$ of the epigraph
of $F$ and consider as $D$ the projection of $E$ onto $X$.
\endproof
 
\lemma{4.2}{Let $d\in\n$ and $X$ be a subset of $\r^d$. Let $H$ be a
Lipschitz multivalued function {}from $X$ to $\Rm$ with non-empty,
compact
and convex values. Then, there exists a sequence $(h_j)$ of Lipschitz
functions {}from $X$ to $\Rm$ such that
$$
H(x)\Eq {\rm cl}\{h_j(x): j\in\n\}
$$
for every $x\in X$, where {\rm cl} denotes the closure in $\Rm$.}
 
\proof  Let $(\xi_j)$ be a dense sequence in $\Rm$; for every
$x\in X$, define
$h_j(x)= P(\xi_j,H(x))$ $\in H(x)$,
where $P$ is the projection map given in Theorem~1.3. Since $P$ and
$H$ are
both lip\-schitzian, so is $h_j$.
\par
Given $\xi\in H(x)$ and $\e > 0$ there exists $\xi_j\in\Rm$ such that
$|\xi
- \xi_j|<\e$. If $L$ denotes a Lipschitz constant for $P$, then
$$
|\xi - h_j(x)|\Eq |P(\xi,H(x)) - P(\xi_j,H(x))|\Le L|\xi - \xi_j|\L L\e\,;
$$
we conclude that
$\xi\in {\rm cl}\{h_j(x):j\in\n\}$.
\endproof
\bigskip
 
We shall now state a result due to G. Bouchitt\'e and M. Valadier
concerning the commutativity property for the operations of integration
and
infimum. To this aim we need the following notion of  $C^1\!$-convexity
which is essentially the notion of $C^1\!$-stability introduced in~[7].
 
\defin{4.3}{Given a positive Radon measure $\lambda$ on $\om$
and a set ${\cal H}$ of $\lambda$-measurable functions {}from $\om$
into
$\Rm$, we say that ${\cal H}$ is $C^1\!$-{\it convex\/} if for every finite
family $(u_i)_{i\in I}$ of elements of ${\cal H}$ and for every
family $(\alpha_i)_{i\in I}$ of non-negative functions of $C^1(\om)\cap
W^{1,\infty}(\om)$ such that $\sum_i\alpha_i= 1$ in $\om$, we have
that
$\sum_i\alpha_i u_i$ belongs to ${\cal H}$.}
\bigskip
 
Let $\lambda$ be a positive Radon measure on $\om$ and ${\cal H}$ be
a
family of
$\lambda$-measurable functions {}from $\om$ to $\Rm$. Then, there
exists a
closed valued $\lambda$-measurable multifunction
$\Gamma\colon\om\rightarrow\Rm$ (i.e., such that
$\Gamma^{-1}(C)=\{x\in\om:\Gamma(x)\cap C\Ne\emptyset\}$ is
$\lambda\!$-measurable for every closed subset $C$ of $\Rm$) with the
following properties (see [29], Proposition~14):
\smallskip
\item{(i)} for every $w\in{\cal H}$ we have $w(x)\in\Gamma(x)$ for
$\lambda$-a.e.\ $x\in\om$;
\smallskip
\item{(ii)} if $\Phi\colon\om\rightarrow\Rm$ is a
closed valued $\lambda$-measurable multifunction such that for every
$w\in{\cal H}$,
$w(x)\in\Phi (x)$ for $\lambda$-a.e.\ $x\in\om$, then $\Gamma
(x)\subseteq
\Phi (x)$ for $\lambda$-a.e.\ $x\in\om$.
\smallskip
This multifunction $\Gamma$
is unique up to $\lambda\!$-equivalence and will be denoted by
\break
$\lambda\hyphen{\displaystyle\essup_{w\in{\cal H}}}\,\{w(\cdot)\}$.
 
\bigskip
The next theorem is taken {}from [7], Theorem~1.
 
\th{4.4}{Let $\lambda$ be a positive Radon measure on
$\om$ and let ${\cal H}$ be a $C^1\!$-convex family of
$\lambda\!$-measurable functions {}from $\om$ into $\Rm$. Let
$f\colon\om\times\Rm\rightarrow {]-\infty,+\infty]}$ be a ${\cal
B}_\lambda(\om)\otimes{\cal B}(\Rm)$-measurable function such that
$f(x,\cdot)$ is convex on $\Rm$ for $\lambda$-a.e.\ $x\in\om$.
Suppose that
$f(\cdot,u(\cdot))\in L^1(\om,\lambda)$ for every $u\in {\cal H}$ and
let
$\Gamma(x)= \lambda\hyphen{\displaystyle\essup_{u\in{\cal
H}}}\{u(x)\}$. Then
$$
\Inf_{u\in{\cal H}}\into f(x,u(x))\,d\lambda\Eq \into
\Inf_{z\in\Gamma(x)} f(x,z)\,d\lambda\,.
$$
}
 
The following technical
result, proven in [14], Lemma~4.2, will be crucial
in the next section.
 
\lemma{4.5}{Let $(w_k)$ be a sequence of functions in $\WP\cap\LL$
converging in $\LL$ to a function $w\in \WP\cap\LL$. Then there exists
a
sequence $(v_k)$ in $\WP\cap\LL$ such that
$v_k(x)\in {\rm co}\{w_1(x),\ldots,w_k(x)\}$ for a.e.\ $x\in\om$ and
$(v_k)$ converges to $w$ strongly in $\WP$.}
 
\lemma{4.6}{Let $\lambda$ be a positive Borel measure on $\om$. Let
$(\g_h)$ and $\g$ be non-negative functions in $L^1(\om,\lambda)$
satisfying
$$
\leqalignno{\g (x)\le\li \g_h &(x)\qquad\hbox{for $\lambda\!$-a.e.\
$x\in
\om\,,$}&(4.1)\cr
\int_\om\g \,d\lambda\geq&\Limsup_{h\rightarrow \infty}
\int_\om\g_h\,
d\lambda\,.&(4.2)\cr}
$$
Then, $(\g_h)$ converges to $\g$ strongly in $L^1(\om,\lambda)$.
}
 
\proof Let us note that, by the Fatou Lemma, (4.2) ensures that
$$
\int_\om\g \,d\lambda\Eq\l \int_\om\g_h \,d\lambda\,.\leqno(4.3)
$$
In view of (4.1)
we have
$$
\g\Le\li (\g_h\wedge\g)\Le
\Limsup_{h\rightarrow \infty}(\g_h\wedge\g)\Le\g\qquad\hbox{on
$\om\,.$}
$$
Thus, the dominated convergence theorem guarantees that
$(\g_h\wedge\g)$
converges to $\g$ in $L^1(\om,\lambda)$, and, in particular
$$
\int_\om\g_h\wedge\g \,d\lambda\ \ \rightarrow\ \ \int_\om\g
\,d\lambda\,.\leqno(4.4)
$$
By noticing that $\g_h + \g =
(\g_h\wedge\g) + (\g_h\vee\g)$, (4.3) and (4.4) permit to conclude that
$$
\int_\om \g_h\vee\g \,d\lambda\ \ \rightarrow\ \ \int_\om\g
\,d\lambda\,;
$$
hence
$(\g_h\vee\g)$ converges to $\g$ in $L^1(\om,\lambda)$,
being $\g_h\vee\g \geq \g$. Now, the conclusion can be obtained by
using again the relation $\g_h = (\g_h\wedge\g) + (\g_h\vee\g) - \g$ on
$\om$.
\endproof
\bigskip
 
We conclude this section with a Dini-type lemma for which we refer,
e.g.,
to~[14], Lemma~4.3.
 
\lemma{4.7}{Let $E$ be a compact subset of $\r^n$ and $(H_k)$
be an increasing (with respect to inclusion) sequence of lower
semicontinuous
multifunctions {}from $E$ to $\Rm$ with closed values. Let $u\in
C^0(E,\Rm)$
such that
$$
u(x)\in{\rm cl}(\bigcup_{k=1}^{\infty} H_k(x))\qquad\hbox{for
every $x\in E\,.$}
$$
Then, for every $r > 0$ there exists $h\in\n$ such that
$B_r(u(x))\cap H_k(x)\Ne\emptyset$ for every $k\geq h$ and for every
$x\in E$.}
\bigskip
 
\parag{5}{Integral representation on ${\bf\WPL}$}
 
The main result of this section is the integral representation of the
functionals of the class $\F$ on the bounded functions of $\WP$
(Theorem~5.4).
\par
Given $F\in\F$, let us introduce the least closed valued multifunction
having the elements of $\dF$ among its selections.
 
\prop{5.1}{Let $F\in\F$ and let $A$ be an open subset of $\om$ with
${\d F(\cdot,A)\Ne\emptyset}$.
Then there exists a closed valued multifunction $K_A$ {}from $A$ to
$\Rm$, unique up to sets of capacity zero, such that
\smallskip
\item{(i)} for every $u\in{\d F(\cdot,A)}$ we have $u(x)\in
K_A(x)$ for q.e.\ $x\in A$;
\smallskip
\item{(ii)} if $H$ is a closed valued multifunction {}from $A$ to $\Rm$
such that
for every $u\in{\d F(\cdot,A)}$ we have $u(x)\in H(x)$ for q.e.\ $x\in
A$,
then
$K_A(x)\subseteq H(x)$ for q.e.\ $x\in A$.
\par\noindent
Moreover, $K_A$ satisfies the following properties:
\smallskip
\item {(iii)} $K_A$ is quasi lower
semicontinuous and $K_A(x)$ is convex for q.e.\ $x\in A$;
\smallskip
\item{(iv)} if $(u_k)$ is a countable dense subset of $\d F(\cdot,A)$,
then
\par\noindent
$$
K_A(x)\Eq {\rm cl}\{u_k(x):k\in\n\}
\Eq {\rm cl}(\bigcup_{k=1}^{\infty}
C_k(x))\qquad\hbox{for q.e.\ $x\in A\,,$}
$$
\item{}where $C_k(x)= \cou$.\smallskip}
 
\proof  The same argument applied in Proposition~3.3 in~[14] works
now   for the subset $\{u|_A:u\in \d F(\cdot,A)\}$ of $W^{1,p}(A,\r^m)$.
\endproof
 
\rem{5.2}{Let $A$ and $A^\p$ be open subsets of $\om$, with $\d
F(\cdot,A)\Ne \emptyset$  and $\d F(\cdot,A^\p)\Ne\emptyset$. If
$K_A$ and
$K_{A^\p}$ are the multifunctions given by the previous proposition,
then
$K_A= K_{A^\p}$ q.e.\ on $A\cap A^\p$.
\par
It is enough to give the proof in the case $A^\p\subseteq A$. Since
$\d F(\cdot,A)\subseteq \d F(\cdot,A^\p)$, the inclusion
$K_A(x)\subseteq
K_{A^\p}(x)$ for q.e.\ $x\in A^\p$ follows immediately {}from property
(i)
satisfied by $K_{A^\p}$ and property (ii) applied to $K_A$ and
$$
H(x)\Eq
\cases{
K_{A^\p}(x)\,,& if $x\in A^\p\,,$
\cr
\Rm\,, & if $x\in A\setminus A^\p\,.$
\cr}
$$
To get the opposite inclusion let us choose $u_0\in\d F(\cdot,A)$. Fix
now
$u\in\d F(\cdot,A^\p)$ and $A^{\pp}\in\A$ with
$A^{\pp}\subset\subset A^\p$.
If $\f$ is a function in $C^1_0(A^\p)$, with $\f= 1$ on $A^{\pp}$ and
$0\leq \f \leq 1$, by the $C^1\!$-convexity and the locality property of
$F$ on
open sets, we have
$$
\eqalign{F(\f u+(1-\f )u_0,A)&\Le F(\f u+(1-\f)u_0,A^\p) +
	F(\f u+(1-\f)u_0,A\setminus {\rm supp}\f)\cr
&\Le F(u,A^\p) + F(u_0,A^\p) + F(u_0,A\setminus {\rm supp}\f)\L
+\infty\,.\cr}
$$
Therefore, $\f u+(1-\f)u_0\in\d F(\cdot,A)$ so that $u(x)\in K_A(x)$ for
q.e.\ $x\in A^{\pp}$. By the arbitrariness of $A^{\pp}$ we deduce that
$u(x)\in
K_A(x)$ for q.e.\ $x\in A^\p$. By applying property (ii) we conclude that
$K_{A^\p}(x)\subseteq K_A(x)$ for q.e.\ $x\in A^\p$.}
 
\lemma{5.3}{Let $s>0$ and let $T_s$ be the orthogonal projection onto
the
ball ${\overline B}_s(0)$ defined in Lemma~2.5. Then for every $u$,
$v\in\WP$
and $A\in\A$
$$
F(u+T_s\circ (v-u),A)\Le F(u,A) + F(v,A)\,.
$$
}
 
\proof It is enough to consider the case $A\subset\subset\om$.
Let $\f\in C_0^1(\om)$, $\f=1$ on $A$, $0\leq \f\leq 1$. By Lemma~3.6
in~[14]
there exists a sequence $(\psi_h)$ of functions in $C^\infty(\r^n)$ such
that
$0\leq \psi_h\leq 1$ and $(\psi_h\f(v-u))$ converges to $T_s\circ[\f(v-
u)]$
weakly in $\WP$ as $h$ goes to $\infty$. Since $F(\cdot,A)$ is weakly
lower
semicontinuous on $\WP$ (recall that $F(\cdot,A)$ is convex) we get
$$
\displaylines{F(u+T_s\circ [\f (v-u)],A) \Le \li F(u+\psi_h\f(v-u),A)\cr
\Le
F(u,A) + F(v,A)\,,\cr}
$$
where in the last estimate we have used the $C^1\!$-convexity of $F$.
Now the conclusion can be obtained by applying the locality property of
$F$ on
open sets.
\endproof
 
\th{5.4}{Let $F\in\F$ with $\dF\cap\LL\Ne\emptyset$. Let
$\nu_0$ be the positive finite Borel measure introduced in
Proposition~2.7 and
let $K= K_{\om}$ be the closed valued multifunction {}from $\om$ to
$\Rm$
given  by
Proposition~5.1 for $A=\om$.  Then, there exist a positive finite Borel
measure
$\mu$ on $\om$, absolutely continuous with respect to capacity, and a
Borel
function $f\colon\om\times\Rm\rightarrow [0,+\infty]$ with the
following
properties:
\smallskip
\item {(i)} for every $x\in\om$, the function $f(x,\cdot)$ is
convex and lower semicontinuous on $\Rm$;
\smallskip
\item {(ii)} for every $u\in\WPL$ and for every $A\in{\cal A}(\om)$
$$
F(u,A)\Eq
\cases{\int_A{f(x,u(x))\,d\mu} + \nu_0(A)\,,& if
                   \ $u(x)\in K(x)$\quad\hbox{for q.e.\ $x\in A\,,$}
\cr
\Extra
\cr
+\infty\,,& otherwise$\,.$
\cr}
$$
\smallskip}
 
\proof
{\it Step 1.} Let $(u_i)$ be a sequence of functions in
$\dF\cap\LL$ which will be specified in Step~2. We construct now the
measure
$\mu$ and the integrand $f$ (see (5.3)) satisfying (i) and we
prove
that for every $A\in{\cal A}(\om)$ and for every $k\in\n$ we have
$$
F(u,A)\Eq\int_Af(x,u(x))\,d\mu+\nu_0(A)\leqno(5.1)
$$
whenever $u\in\WPL$ and $u(x)\in\cou$ for
a.e.\ $x\in A$.
\par
Fix quasi continuous Borel measurable
representatives of $(u_i)$. For every $x\in\om$ define
$C_k(x)=\cou$. By Theorem~3.7, for every
$k\in\n$ there exist a positive finite Borel measure $\mu_k$ on $\om$,
absolutely
continuous with respect to capacity, and a function
$f_k\colon\om\times\Rm\rightarrow [0,+\infty]$ such that
\smallskip
\item {(a)} for every $x\in\om$ the function $f_k(x,\cdot)$
is convex and lower semicontinuous on $\Rm$;
moreover, the restriction of
$f_k(x,\cdot)$ to $C_k(x)$ is continuous for $\mu_k\!$-a.e.\ $x\in\om$;
\smallskip
\item {(b)} $f_k$ is
${\wh{\cal B}}(\om)\otimes{\cal B}(\Rm)$-measurable;
\smallskip
\item {(c)}
for every $u\in\WPL$ such that $u(x)\in C_k(x)$ for a.e.\ $x\in\om$ the
function $f_k(\cdot,u(\cdot))$ is $\mu\!$-measurable on $\om$ and for
every
$A\in\A$
\par\vskip-10pt
$$
F(u,A)\Eq\int_A{f_k(x,u(x))\,d\mu_k} + \nu_0(A)\,.\leqno(5.2)
$$
By a standard cut-off argument we obtain that (5.2) still holds if
$u(x)\in
C_k(x)$ for  a.e.\ $x\in A$.
\par
Let $\mu$ be a positive finite Borel measure on $\om$ absolutely
continuous
with respect to capacity and such that $\mu_k\ll\mu$ for every
$k\in\n$
(for instance, take $\mu(B)=\sum_{k=1}^{\infty}
2^{-k}{\mu_k(B)\over{\mu_k(\om)}}$ for every $B\in\B)$. Define
$$
g_k(x,\xi)\Eq f_k(x,\xi){d\mu_k\over{d\mu}}(x)\,,
$$
where $d\mu_k/d\mu$ is a fixed ($\mu\!$-measurable) representative
of the
Radon-Nikodym derivative of $\mu_k$ with respect to $\mu$.
By Proposition~3.8, there exists a set $N\in\B$ with $\mu(N)=0$ such
that
$g_k(x,\xi)=g_{k+1}(x,\xi)$ for every $k\in\n$, $x\in\om\setminus N$
and
$\xi\in C_k(x)$. Hence, we can define $g\colon
\om\times\Rm\rightarrow
[0,+\infty]$ as
$$
g(x,\xi)\Eq
\cases{
g_k(x,\xi)\,, & if\ \ $x\in\om\setminus N$ and $\xi\in
C_k(x) $ for some $k\in\n\,,$\cr
+\infty\,, & otherwise.
\cr}
$$
Since $C_k(x)=\{\sum_{i=1}^k\lambda^iu_i(x): \lambda\in\Sigma_k\}$,
by
Theorem~III.9 and Proposition~III.13 in~[10], the graph of $C_k$
belongs
to $\B\otimes{\cal B}(\Rm)$. Recalling the definition of $g_k$ it follows
that
$g$ is ${\wh{\cal B}}(\om)\otimes{\cal B}(\Rm)$-measurable. An easy
check gives
the convexity of $g(x,\cdot)$ on $\Rm$ for every $x\in\om$.
\par
Now, for every $x\in\om$ let us set $h(x,\cdot)= {\rm sc}^-g(x,\cdot)$,
where ${\rm sc}^-g(x,\cdot)$ denotes the lower semicontinuous
envelope of
$g(x,\cdot)$. It turns out that
$$
h(x,\xi)=\sup_{s\in\n}g_s(x,\xi)\,,
$$
where $g_s(x,\xi)= {\displaystyle \inf_{\eta\in\Rm}[g(x,\eta)+ s|\xi
- \eta|]}$. By Remark~3.6, for every $s\in\n$ there exists a set
$Z_s\in\B$, with $\mu (Z_s)= 0$, such that $g_s|_{(\om\setminus
Z_s)\times\Rm}$ is a Borel function. Set $Z=\bigcup_{s\in\n}Z_s$; then
$\mu(Z)=0$ and $h|_{(\om \setminus Z)\times\Rm}$ is Borel
measurable. Now we
are in a position to define the function $f$ as
$$
f(x,\cdot)\Eq
\cases{
h(x,\cdot)\Eq {\rm sc}^-g(x,\cdot)\quad &if\ \ $x\in\om\setminus Z$
\cr
0&if\ \ $x\in Z\,.$
\cr}\leqno (5.3)
$$
Then, $f$ is a Borel function on $\om\times\Rm$ and satisfies (i) (see
[28], Theorem~7.4).
\par
Let us prove that for $\mu\!$-a.e.\ $x\in\om$
$$
f(x,\cdot)\Eq g_k(x,\cdot)\qquad\hbox{on $C_k(x)$.}\leqno(5.4)
$$
Let us fix $x\in\om\setminus(N\cup Z)$ and $k\in\n$. Let $H(x)$ be the
affine
hull of $\bigcup_{k=1}^\infty C_k(x)$. As the sequence $(C_k(x))$ is
increasing, there exists $h\Ge k$ such that the interior of $C_h(x)$
relative
to $H(x)$ is non-empty. Since $g(x,\cdot)=g_h(x,\cdot)$ on $C_h(x)$, and
the
restriction of $g_h(x,\cdot)$ to $C_h(x)$ is continuous, we have
$f(x,\cdot)=g_h(x,\cdot)$ on the interior of $C_h(x)$ relative to $H(x)$.
As
$C_h(x)$ is a polytope, the restriction of $f(x,\cdot)$ to $C_h(x)$ is
continuous (see [28], Theorem~10.2), hence $f(x,\cdot)=g_h(x,\cdot)$ on
$C_h(x)$. Since $C_k(x)\subseteq C_h(x)$ and
$g_k(x,\cdot)=g_h(x,\cdot)$ on
$C_k(x)$, we conclude that (5.4) is satisfied.
\par
Let us now prove (5.1). Fix $k\in\n$, $A\in{\cal A}(\om)$ and
$u\in\WPL$ with $u(x)\in C_k(x)$ for a.e.\ $x\in A$. For every
$0<\sigma< 1$
let us define $u_\sigma= \sigma u+(1-\sigma)u_0$, with $u_0= {1\over
k}\sum_{i=1}^k u_i$. Then $u_\sigma (x)\in$ ri$C_k(x)$ for q.e.\ $x\in
A$.
Therefore, by (5.4) and the integral representation formula (5.2)
satisfied by
$f_k$, we get %
$$
F(u_\sigma,A)\Eq \int_A f(x,u_\sigma(x))\,d\mu\ +\nu_0(A)\,.
$$
Since every lower semicontinuous proper convex function is continuous
along
line segments (see [28], Corollary 7.5.1) it turns out that
$$
\eqalign{\lim_{\sigma\rightarrow 1^-} F(u_\sigma,A)&\Eq F(u,A) \cr
\lim_{\sigma\rightarrow 1^-} \int_A f(x,u_\sigma(x))\,d\mu &\Eq
\int_A f(x,u(x))\,d\mu\qquad\hbox{for every $x\in\om\,.$}
\cr}
$$
We thus obtain (5.1).
\bigskip
\noindent
{\it Step 2.} We choose now a suitable sequence $(u_i)$, dense in
$\dF\cap\LL$,
to which Step~1 will be applied.
\par
Let ${\cal D}$ be a countable base for the open subsets of
$\om$, closed under finite unions. For every $A\in{\cal D}$ we can
apply
Lemma~4.1 to
$F(\cdot,A)$ on  $\dF\cap\LL$ (with the $\WP$ topology); this yields
the existence of a set ${\cal G}_A\subseteq\dF\cap\LL$
such that for every $u\in\dF\cap\LL$ there exists a
sequence $(u_h)$ in ${\cal G}_A$ satisfying
$$
\displaylines{u_h \ \rightarrow\
u\qquad\hbox{strongly in $\WP$,}\cr
F(u_h,A)\ \rightarrow\  F(u,A)\qquad\hbox{in
$\r\,$.}\cr}
$$
Let $(u_i)$ be an enumeration of $\bigcup_{A\in{\cal D}}{\cal
G}_A$; starting {}from $(u_i)$ we then construct by means of Step~1 a
Borel function $f\colon\om\times\Rm\rightarrow [0,+\infty]$ satisfying
(i) and
(5.1).
\bigskip
\noindent
{\it Step 3.} Let us prove that for every $u\in\dF\cap\LL$ and
for every $A\in{\cal A}(\om)$
$$
F(u,A)\Ge\int_Af(x,u(x))\,d\mu + \nu_0(A)\,.\leqno(5.5)
$$
\par
Fix $u\in\dF\cap\LL$ and $A\in{\cal D}$. By Step~2 it is
possible to extract a sequence $(u_{i_{\scriptscriptstyle h}})$ {}from
$\{u_i:
i\in\n\}$ such that
$$
\displaylines{{u_{i_h}}\ \rightarrow\
u\qquad\hbox{q.e.\ in $\om$}\qquad\hbox{(hence $\mu$-a.e.)$\,,$}\cr
F(u_{i_h},A)\ \rightarrow\ F(u,A)\qquad\hbox{in $\r$}\,.\cr}
$$
Therefore, by (5.1)
$$
F(u,A)\Eq \l\int_A f(x,u_{i_{\scriptscriptstyle h}}(x))\,d\mu +
\nu_0(A)\,;
$$
by the Fatou Lemma and the lower semicontinuity of $f(x,\cdot)$ we
get (5.5) for every $A\in {\cal D}$. The result for an arbitrary $A\in\A$
can
be obtained by approximation, since each $A\in\A$ is the union of an
increasing sequence of elements of ${\cal D}$ (recall that ${\cal D}$ is
closed under finite unions).
\bigskip
\noindent
{\it Step 4.} It is now easy to prove that for every $u\in\WPL$ and
$A\in{\cal A}(\om)$ the inequality (5.5) holds.
\par
Given $A\in{\cal A}(\om)$ and $u\in\d F(\cdot,A)\cap\LL$,
let $G\subset\subset A$ and $\f\in C^1_0(A)$ with
$\f= 1$ on $G$ and $0\leq \f\leq 1$. Set $u_\f = \f u + (1-\f)w$, where
$w$ belongs to $\dF\cap\LL$, which is non-empty by assumption. By
the convexity and the locality property of $F$ we
have $F(u_\f,\om)\L +\infty$.
Therefore, Step~3 applied to $u_\f$ yields
$$
F(u_\f,G)\Ge\int_G f(x,u_\f(x))\,d\mu + \nu_0(G)\,;
$$
since $\f = 1$ on $G$ we get
$$
F(u,G)\Ge\int_G f(x,u(x))\,d\mu + \nu_0(G)\,.
$$
As $G\subset\subset A$ is arbitrary, the conclusion is easily achieved.
\bigskip
\noindent
{\it Step 5.} Let $K=K_\om$ be the closed valued multifunction {}from
$\om$ to
$\Rm$ given by Proposition~5.1 for $A=\om$. The aim is now to prove
that for
every $A\in{\cal A}(\om)$ and $u\in\WPL$ such that $u(x)\in K(x)$ for
q.e.\
$x\in A$, we have %
$$
F(u,A)\Le\int_Af(x,u(x))\,d\mu + \nu_0(A)\,.\leqno(5.6)
$$
\par
Recall that for every $k\in\n$ and $x\in\om$ we have $C_k(x)=
{\rm co}\{u_1(x),\ldots,u_k(x)\}$, where $(u_i)$ is the
sequence given in Step~2. By
Lemmas~2.5 and 5.3, $\dF\cap L^\infty(\om,\r^m)$, which is non-empty
by
assumption, is dense in $\dF$. Hence, $(u_i)$ is dense in $\dF$ and, by
Proposition~5.1, $K(x)= {\rm cl}\bigl(\bigcup_{k=1}^{\infty}
C_k(x)\bigr)$ for
q.e.\ $x\in\om$.
\par
Fix $u\in\WPL$ and $A\in{\cal A}(\om)$ such that $u(x)\in K(x)$ for
q.e.\
$x\in A$. Clearly, we can assume that the right-hand side in (5.6) is
finite. Moreover, we can consider open sets $A\subset\subset\om$ with
smooth boundary, so that there exists an extension operator
$W^{1,p}(A)\rightarrow W^{1,p}(\om)$.
 
\bigskip
In a first moment we work with the additional assumption that $u(x)$ is
in the closure of $\bigcup_{k=1}^{\infty} C_k(x)$ ``uniformly" for $x\in
A$;
more precisely, given a sequence $(r_h)$ of positive numbers decreasing
to $0$, we require that for every $h\in\n$ there exists $n_h\in\n$ such
that
$$
B_{r_h/2}(u(x))\cap C_k(x)\Ne\emptyset
\quad\hbox{for q.e.\
$x\in A$}\leqno(5.7)
$$
for every $k\Ge n_h$.
\par
To achieve (5.6) we look for a sequence $(v_h)$ of functions in
$W^{1,p}(A,\r^m)\cap L^\infty(A,\r^m)$ and a strictly increasing
sequence $(k_h)$ of positive integers such
that $v_h\in C_{k_h}(x)$ for q.e.\ $x\in A$, $(v_h)$ converges to $u$ in
$W^{1,p}(A,\r^m)$ and
$$
\displaystyle\limsup_{h\rightarrow \infty}\int_A
f(x,v_h(x))\,d\mu\Le\int_A
f(x,u(x))\,d\mu\,.\leqno(5.8)
$$
Indeed, as $\partial A$ is smooth, we can assume that $(v_h)$ is a
sequence in $\WPL$ converging in $\WP$ to a function
$v$ such that $v=u$ a.e.\ on $A$.
By the lower semicontinuity of
$F(\cdot,A)$ and the integral representation (5.1) obtained in Step~1, we
can
then conclude
$$
\eqalign{F(u,A)\Eq F(v,A)&\Le
\displaystyle\liminf_{h\rightarrow \infty}F(v_h,A)\cr &\Le
\displaystyle\Limsup_{h\rightarrow \infty}\bigl(\int_A
f(x,v_h(x))\,d\mu +
\nu_0(A)\bigr) \Le\int_A f(x,u(x))\,d\mu + \nu_0(A)\,.\cr}
$$
\par
Let us first construct a sequence $(w_h)$ of functions in
$W^{1,p}(A,\r^m)\cap L^\infty(A,\r^m)$ and a strictly increasing
sequence $(k_h)$ of positive integers,
with the following properties:
$$
\eqalign{w_h(x)\in C_{k_h}(x)\qquad\hbox{for q.e.\ $x\in A$}
\cr
w_h\ \rightarrow\ u\qquad\hbox{uniformly on $A$}
\cr
\displaystyle\Limsup_{h\rightarrow \infty}\int_A f(x,w_h(x))\,d\mu
\Le \int_A f(x,u(x))\,d\mu\,.
\cr}
\leqno(5.9)
$$
To this aim let us prove that for every $h\in\n$ there exists $k_h\in\n$
such
that
$$
\displaystyle\inf_{w\in{\cal H}^h_k}\int_A f(x,w(x))\,d\mu\L \int_A
f(x,u(x))\,d\mu + r_h\,,\leqno(5.10)
$$
for every $k\Ge k_h$, where
$$
{\cal H}^h_k\Eq\{w\in W^{1,p}(A,\r^m)\cap L^\infty(A,\r^m):
w(x)\in{\overline
B}_{r\scriptstyle_h}(u(x))\cap C_k(x)\ \,\hbox{for q.e.\ $x\in
A$}\}\,.
$$
Let us fix $h\in\n$ and let $n_h$ be as in (5.7). For every fixed $k\ge
n_h$
we want to apply Theorem~4.4 to the set ${\cal H}^h_k$. For this
purpose let
us verify that
$$
\mu\hyphen{\displaystyle\essup_{w\in{\cal H}^h_k}}\{w(x)\}\Eq
{\overline
B}_{r\scriptstyle_h}(u(x))\cap C_k(x)\qquad\hbox{for q.e.\ $x\in
A\,.$}\leqno(5.11)
$$
Let $\Xi =\{(\xi,\vxi)\in(\Rm)^{k+1}:d(\xi,\cox)\Le r_h/2\}$
and let $H$ be the multivalued function {}from $\Xi$ to $\Rm$ defined
by
$$
H(\xi,\vxi)\Eq {\overline B}_{r\scriptstyle_h}(\xi)\cap\cox\,.
$$
By Theorem~1 in~[24], $H$ is lipschitzian. Hence, we can apply
Lemma~4.2 to
$H$ obtaining a sequence $(h_j)$ of Lipschitz functions {}from $\Xi$ to
$\Rm$
such that
$$
H(\xi,\vxi)\Eq {\rm cl}\displaystyle\bigl(\bigcup_{j=1}^{\infty}
\{h_j(\xi,\vxi)\}\bigr)\,.
$$
Since $B_{r_h/2}(u(x))\cap C_k(x)\neq\emptyset$ for q.e.\
$x\in A$, we can define $z_j=h_j(u,u_1,\cdots,u_k)$\ q.e.\ on $A$
for every $j\in\n$. By Remark~1.6., $z_j\in W^{1,p}(A,\r^m)\cap
L^\infty(A,\r^m)$. Thus $z_j\in{\cal H}^h_k$, and
$$
{\overline B}_{r\scriptstyle_h}(u(x))\cap C_k(x)\Eq H(u(x),
u_1(x),\ldots,u_k(x))\Eq {\rm cl}\bigl(\bigcup_{j=1}^{\infty}\{
z_j(x)\}\bigr)
$$
for q.e.\ $x\in A$. Hence, (5.11) holds.
\par\noindent
Moreover, since every $w\in{\cal H}^h_k$ is a convex combination of
$\uxi$, we
have
$$
\int_A f(x,w(x))\,d\mu \Le \int_A \displaystyle\sum_{i=1}^k
f(x,u_i(x))\,d\mu
\L +\infty\,.
$$
We can now apply Theorem~4.4; by (5.11), for every $k\ge n_h$ we
have
$$
\inf_{w\in{\cal H}^h_k}\int_A f(x,w(x))\,d\mu \Eq \int_A\
\inf_{\xi\in C_k^h(x)}f(x,\xi)\,d\mu\,,\leqno(5.12)
$$
where $C_k^h(x)={\overline B}_{r\scriptstyle_h}(u(x))\cap C_k(x)$. Since
$u(x)\in K(x)$ for q.e.\ $x\in A$, in view of the continuity property
along
line segments for a proper, lower semicontinuous convex function, for
q.e.\
$x\in A$ we can approximate $u(x)$ by a sequence $(\xi_k(x))$ in
${\rm ri}K(x)$
such that
$$
f(x,u(x))\Eq \lim_{k\rightarrow \infty}f(x,\xi_k(x))\,.
$$
As ${\rm ri}K(x)\subseteq\bigcup_{k=1}^{\infty}C_k(x)$ (see [28],
Theorem~6.3), we can suppose that $\xi_k(x)\in C_k(x)$ for every
$k\in\n$.
Thus, for every $h\in\n$ we have %
$$
\inf_{k\in\n}\ \inf_{\xi\in C_k^h(x)} f(x,\xi)\Le
f(x,u(x))\qquad\hbox{for q.e.\ $x\in A$}\,.\leqno(5.13)
$$
Moreover, for every
$k\geq n_h$ the set $C_k^h(x)$ is
non-empty for q.e.\ $x\in A$. Therefore,
the convexity of $f$ ensures that
$$
\int_A {\inf_{\xi\in C_k^h(x)} f(x,\xi)\,d\mu }\Le
\int_A {\displaystyle\sum_{i=1}^k f(x,u_i(x))\,d\mu} \L +\infty\,;
$$
by (5.13) and by the monotone convergence theorem it follows that
$$
\inf_{k\in\n}\int_A {\inf_{\xi\in C_k^h(x)} f(x,\xi)\,d\mu} \Le
\int_A {f(x,u(x))\,d\mu}\,.
$$
This inequality, together with (5.12), proves (5.10).
\par
Let $(k_h)$ be the sequence given in (5.10) which we can assume to be
strictly
increasing. For every $h\in\n$, by (5.10) there exists a
function $w_h\in W^{1,p}(A,\r^m)\cap L^\infty (A,\r^m)$ such that
$$
\eqalign{& w_h(x)\in{\overline B}_{r\scriptstyle_h}(u(x))\cap
C_{k_h}(x)\qquad\hbox{for q.e.\ $x\in A\,,$}\cr
&\int_A f(x,w_h(x))\,d\mu \Le \int_A f(x,u(x))\,d\mu +
r_h\,.\cr}
$$
It is easy to verify that $(w_h)$ satisfies the properties in (5.9).
\bigskip
 
Let us set now $\g_h=f(\cdot,w_h(\cdot))$ and $\g =f(\cdot,u(\cdot))$.
We claim that
$$
\g_h\rightarrow\g\qquad\hbox{strongly in
$L^1(A,\mu)\,.$}\leqno(5.14)
$$
Indeed, as $(w_h)$ converges to $u$ q.e.\ on $A$, by the lower
semicontinuity
of $f(x,\cdot)$ we get $\g (x)\le{\displaystyle\liminf_{h\rightarrow
\infty}\g_h(x)}$. By (5.9) and Lemma~4.6, it follows that $(\g_h)$
converges to
$\g$ in the strong topology of $L^1(A,\mu)$.
\par
In view of (5.14) it is not restrictive to assume that for every $h\in\n$
$$
\int_A |\g_h - \g|\,d\mu\Le {1\over{2^h}}\,.\leqno(5.15)
$$
At this point let us apply Lemma~4.5 to the sequence $(w_j)_{j\ge h}$
for every $h\in\n$. We obtain a sequence $(v_{h,j})_{j\ge h}$ of
functions
in $W^{1,p}(A,\Rm)\cap L^\infty(A,\Rm)$ such that
$$
\displaylines{v_{h,j}(x)\in{\rm co}\{w_h(x),w_{h+1}(x),
\ldots,w_j(x)\}\qquad\hbox{for q.e.\ $x\in A\,,$}\cr
v_{h,j}\rightarrow u\qquad\hbox{strongly in $W^{1,p}(A,\Rm)$ as
$j\rightarrow \infty\,.$}\cr}
$$
By a standard argument we can find a strictly increasing sequence
$(j_h)$ of positive integers such that $(v_{h,j_h})$ converges to $u$
in $W^{1,p}(A,\Rm)$. Define $v_h=v_{h,j_h}$ for every $h\in\n$. Then
$v_h\in
W^{1,p}(A,\Rm)\cap L^\infty(A,\Rm)$ and
$$
\displaylines{v_h(x)\in{\rm co}\{w_h(x),\ldots,
w_{j_h}(x)\}\qquad\hbox{for q.e.\ $x\in A\,,$}\cr
v_h\rightarrow u\qquad\hbox{strongly in $W^{1,p}(A,\Rm)\,$}.\cr}
$$
In particular, a suitable sequence $(k_h)$ exists such that $v_h(x)\in
C_{k_h}(x)$ for q.e.\ $x\in A$. Now we only need to verify that (5.8)
holds
for the sequence $(v_h)$ just obtained. By Lemma~3.2 we can write
$v_h(x)=
\sum_{i=h}^{k_h}\psi^i_h(x)w_i(x)$ for q.e.\ $x\in A$, where
$\psi_h\colon
A\rightarrow \Sigma_{k_h-h+1}$ are $\mu$-measurable. Let us now
make use of
the convexity of $f$, together with (5.15):
$$
\displaylines{\int_A f(x,v_h(x))\,d\mu \Le\sum_{i=h}^{k_h}\int_A
\psi_h^i(x)|\g_i(x) - \g(x)|\,d\mu + \int_A\g(x)\,d\mu\cr
\Le\sum_{i=h}^{k_h}{1\over{2^i}} + \int_A\g(x)\,d\mu\Le
\int_A\g(x)\,d\mu
+ {1\over{2^h}}\,.\cr}
$$
This implies
$$
\Limsup_{h\rightarrow \infty}\int_A f(x,v_h(x))\,d\mu\Le \int_A
f(x,u(x))\,d\mu\,.
$$

\bigskip
Finally, let us remove the additional assumption (5.7).
Fix $G\subset\subset A$ and a sequence $(r_h)$ of positive real
numbers
decreasing to $0$. For every $\e > 0$ there exists an open set
$A_\e\subseteq\om$, with ${\rm cap}(A_\e,\om)<\e$, such that
$u_i|_{\om\setminus A_\e}$ and $u|_{\om\setminus A_\e}$ are
continuous for
every $i\in\n$. In particular, the multifunction $C_k$ is continuous on
$\om\setminus A_\e$ with respect to the Hausdorff metric. By
Lemma~4.7 for
every $h\in\n$ there exists $n_h^\e\in\n$ such that
$B_{r_h/2}(u(x))\cap C_k(x)\Ne
\emptyset$ for every $k\Ge n^\e_h$ and for every $x\in G\setminus
A_\e$.
Let $z_\e$ be the capacitary potential of $A_\e$ and $u_\e=
(1-z_\e)u + z_\e u_1$, where $u_1$ is the first term of the sequence
$(u_i)$.
Then one can easily check that $(u_\e)$ converges to $u$ in $\WP$, that
$u_\e\in K(x)$ for q.e.\ $x\in G$, and that for every $h\in\n$ there
exists
$n^\e_h\in\n$ such that $B_{r_h/2}(u_\e (x))\cap
C_k(x)\Ne \emptyset$ for q.e.\ $x\in G$ and for every $k\Ge n_h^\e$.
Therefore we can apply the previous result for $u_\e$ and $G$ in place
of $u$
and $A$; this gives
$$
F(u_\e,G)\Le\int_G f(x,u_\e (x))\,d\mu + \nu_0(G)\,.
$$
Since
$$
\int_G f(x,u_\e(x))\,d\mu \Le \int_A \bigl[(1-z_\e(x))f(x,u(x))
+ z_\e(x) f(x,u_1(x))\bigr]\,d\mu
$$
and $\int_A f(x,u(x))\,d\mu < +\infty$ by assumption, the lower
semicontinuity of $F$ and the dominated convergence theorem imply
$$
F(u,G)\Le \liminf_{\e\rightarrow 0^+}F(u_\e ,G)\Le \int_A
f(x,u(x))\,d\mu +
\nu_0(A)\,.
$$
Taking the supremum for $G\subset\subset A$ we get (5.6).
\bigskip
\noindent
{\it Step 6.} In view of Step~4 and Step~5 we get
$$
F(u,A)\Eq\int_Af(x,u(x))\,d\mu+\nu_0(A)
$$
for every $u\in\WPL$ with $u(x)\in K(x)$ for q.e.\ $x\in A$. Property
(ii)
now follows by taking into account that if
$u\in\d F(\cdot,A)$ then $u(x)\in K(x)$ for q.e.\ $x\in A$ by
Remark~5.2.
\endproof
\bigskip
 
\parag{6}{Integral representation on ${\bf\WP}$}
 
We now eliminate (Theorem~6.1) the restrictive condition $u\in\LL$
considered in the previous section. Furthermore, Proposition 6.3 will
allow us
to treat in a unified way both cases of the representation formula
established
in Theorem~5.4(ii). Thus, we achieve (Theorem~6.5) the conclusive
integral
representation theorem, which is the main result of the paper.
\bigskip
 
Given $F\in\F$ let us define
$$
\overline{\nu}(B)\Eq\inf\{F(u,B):u\in\WP\}\,,\leqno(6.1)
$$
for every $B\in\B$. It is easily seen that the proof of Proposition 2.7
still works for the set function $\overline{\nu}$ on every
$\om^\p\in\A$ with
$\d F(\cdot,\om^\p)\ne\emptyset$. Therefore, on such sets,
$\overline{\nu}$ is
a  positive finite Borel measure.
 
\th{6.1}{Let $F\in\F$ and assume $\dF\Ne\emptyset$. Then the
conclusions
of Theorem~5.4 still hold with $u$ (in item (ii)) ranging all over $\WP$
and
$\nu_0$ replaced by the measure $\overline{\nu}$ defined in (6.1).
}
 
\proof For every $v\in\WP$ and for every $B\in\B$ let us define
$$
\displaylines{X_v\Eq\{u\in\WP: u-v\in\WPL\}\,,\cr
               \nu_v(B)={\rm inf}\{F(u,B):u\in X_v\}\,.\cr}
$$
By a suitable application of Theorem~5.4, it turns out that for every
$v\in\dF$ there exist a positive finite Borel
measure $\mu_v$ on $\om$, absolutely continuous
with
respect to capacity, and a Borel function
$f_v\colon\om\times\Rm\rightarrow
[0,+\infty]$ such that
\smallskip
\item {(i)} for every $x\in\om$ the function $f_v(x,\cdot)$ is convex
and lower semicontinuous on $\Rm$;
\smallskip
\item {(ii)} for every $u\in X_v$ and for every $A\in{\cal A}(\om)$
\par\noindent
$$
F(u,A)\Eq
\cases{
\int_A{f_v(x,u(x))\,d\mu_v} + \nu_v(A)\,,& if
                   \ $u(x)\in K(x)$\quad for q.e.\ $x\in A\,,$
\cr
\Extra
\cr
+\infty\,,& otherwise$\,,$\cr}\leqno(6.2)
$$
\item {}
where $K= K_\om$ is the closed valued multifunction {}from $\om$ to
$\Rm$
given by Proposition 5.1 for $A=\om$.
 
\smallskip\noindent
{\it Step 1.} Let us show first that for every $v\in\dF$,
$u\in\WP$, and $A\in\A$ we have
$$
F(u,A)\L+\infty\quad\hbox{if and only if}\ \
\cases{
u(x)\in K(x)\quad\hbox{for q.e.\ $x\in A\,,$}
\cr
\Extra
\cr
\int_A{f_v(x,u(x))\,d\mu_v} + \nu_v(A)\L +\infty\,.\cr}\leqno(6.3)
$$
By the definition of $K$ and Remark~5.2, if $F(u,A)<+\infty$ then
$u(x)\in
K(x)$ for q.e.\ $x\in A$. Hence, let us assume that $u(x)\in K(x)$
for q.e.\ $x\in A$ and prove that $F(u,A)<+\infty$ if and only if
$\int_A{f_v(x,u(x))\,d\mu_v} + \nu_v(A)\L +\infty$.
\par
For every $k\in\n$, let $T_k\colon\Rm\rightarrow\Rm$ be the
orthogonal
projection onto the ball ${\overline B}_k(0)$;
by Lemma~2.5, for every
$w\in\WP$ the function $T_k\circ w$ belongs to $\WPL$, and the
sequence
$(T_k\circ w)$
converges to $w$ in the strong topology of $\WP$ as $k$ tends to
$\infty$.
\par
For every $k\in\n$ let us set $u_k=v+T_k\circ (u-v)$. By (6.2)
we have
$$
F(u_k,A)\Eq \int_A{f_v(x,u_k(x))\,d\mu_v} + \nu_v(A)\,.\leqno(6.4)
$$
\par
Assume now $F(u,A)<+\infty$. By (6.4) and Lemma~5.3 we have
$$
\int_A{f_v(x,u_k(x))\,d\mu_v} + \nu_v(A)\Le F(u,A)+F(v,A)\,.
$$
Since, up to a subsequence, $(u_k)$ converges to $u$ q.e.\ on
$\om$, the Fatou Lemma and the lower semicontinuity of $f_v(x,\cdot)$
ensure that
$$
\int_A f_v(x,u(x))\,d\mu_v + \nu_v(A)\Le F(v,A) + F(u,A)\L +\infty\,.
$$
\par Conversely, assume $\int_A f_v(x,u(x))\,d\mu_v + \nu_v(A)\L
+\infty\,.$
For every $k\in\n$, by (6.2), (6.4), and by the convexity of
$f_v(x,\cdot)$ we
get
$$
\displaylines{\llap(6.5)\hfill F(u_k,A)\Le\int_A f_v(x,u(x))\,d\mu_v +
\int_A f_v(x,v(x))\,d\mu_v +
\nu_v(A)\Eq\hfill
\cr
\Eq\int_A f_v(x,u(x))\,d\mu_v + F(v,A)\,.\cr}
$$
Hence, by the lower semicontinuity of $F(\cdot,A)$ we conclude that
$F(u,A)< +\infty$.
\medskip
\noindent
{\it Step 2.} Let us fix $A\in\A$ and $u$, $v\in\dF$. We claim that
$$
F(u,A)\Eq \int_A f_v(x,u(x))\,d\mu_v +
\nu_v(A)\,.\leqno(6.6)
$$
\par
Let us show first that for every $w\in X_u$ with
$w(x)\in {\rm co}\{u(x),v(x)\}$ for q.e.\ $x\in\om$, we have
$$
F(w,A)\Le
\int_A f_v(x,w(x))\,d\mu_v + \nu_v(A)\,.\leqno(6.7)
$$
Let us fix $w\in X_u$ with $w(x)\in {\rm co}\{u(x),v(x)\}$ for q.e.\
$x\in\om$, and let $u_k= v + T_k\circ (w-v)$ for $k\in\n$. By the
lower semicontinuity of $F(\cdot,A)$ and by (6.2)
$$
F(w,A)\Le \liminf_{k\rightarrow \infty}F(u_k,A)\Eq
\liminf_{k\rightarrow \infty}\int_Af_v(x,u_k(x))\,d\mu_v +
\nu_v(A)\,.\leqno(6.8)
$$
Note that $(u_k)$ converges to $w$ q.e.\ on $\om$.
Since $u_k(x)$ is on the segment with endpoints $u(x)$ and $v(x)$,
by the convexity of $f_v$ it turns out that
$f_v(x,u_k(x))\Le f_v(x,v(x)) + f_v(x,u(x))$
for q.e.\  $x\in\om$. {}From (6.3) we have
$\int_A f_v(x,v(x))\,d\mu_v< +\infty$ and\break ${\int_A
f_v(x,u(x))\,d\mu_v<
+\infty}$. Hence, by the continuity property of $f_v(x,\cdot)$
along line segments ([28], Corollary 7.5.1) and the dominated
convergence
theorem
$$
\lim_{k\rightarrow \infty}\int_A f_v(x,u_k(x))\,d\mu_v\Eq \int_A
f_v(x,w(x))\,d\mu_v\,.\leqno(6.9)
$$
In view of (6.8), this implies (6.7).
\par
{}From (6.7) with $w=u$ we obtain
$$
F(u,A)\Le\int_A f_v(x,u(x))\,d\mu_v+\nu_v(A)\,.
$$
Let us now prove the opposite inequality.
\par
If in (6.7) we apply (6.2) to represent $F(w,A)$ we obtain
$$
\int_A f_u(x,w(x))\,d\mu_u + \nu_u(A)\Le
\int_A f_v(x,w(x))\,d\mu_v + \nu_v(A)\,.
$$
By exchanging now the
roles of $u$ and $v$ we obtain that for every $w\in X_v$ with
$w(x)\in {\rm co}\{u(x),v(x)\}$ for q.e.\ $x\in\om$
$$
\int_A f_v(x,w(x))\,d\mu_v + \nu_v(A)\Le \int_A f_u(x,w(x))\,d\mu_u
+
\nu_u(A)\,.
$$
Now, if we take $w= v + T_k\circ (u - v)$ and argue as for
(6.9), by the Fatou Lemma we get
$$
\eqalign{\int_A f_v(x,u(x))\,d\mu_v + \nu_v(A)&\Le
\liminf_{k\rightarrow
\infty}\int_A f_u(x,v(x) + T_k(u(x)-v(x))\,d\mu_u + \nu_u(A)\cr
&\Eq \int_A f_u(x,u(x))\,d\mu_u + \nu_u(A)= F(u,A)\,.\cr}
$$
\medskip
\noindent
{\it Step 3.}
For every $v\in\dF$, $A\in\A$, and $u\in\d F(\cdot,A)$ it turns out
that
$$
F(u,A)\Eq \int_A f_v(x,u(x))\,d\mu_v + \nu_v(A)\,.\leqno(6.10)
$$
This follows by applying the same argument used in Step~4 of
Theorem~5.4.
\medskip
\noindent
{\it Step 4.} Since $\dF\ne\emptyset$, there exists a function $v$ for
which (6.10) holds for every $A\in\A$ and $u\in\d F(\cdot,A)$.
Finally, we obtain that for every $u\in\WP$ and $A\in\A$
$$
F(u,A)\Eq
\cases{
\int_A{f_v(x,u(x))\,d\mu_v} + \nu_v(A)\,,& if
                   \ $u(x)\in K(x)$\quad\hbox{for q.e.\ $x\in A\,,$}
\cr
\Extra
\cr
+\infty\,,& otherwise$\,.$\cr}\leqno(6.11)
$$
Indeed, if $u(x)\in K(x)$ for q.e.\ $x\in A$
and $u\notin \d F(\cdot,A)$, then by (6.3) we have
$\int_A f_v(x,u(x))\,d\mu_v  + \nu_v(A)= +\infty$.
\medskip
\par
So far we have proved the integral representation by means of any of
the
measures $\nu_v$ with $v\in\dF$. We claim that for every $v\in\dF$
and
$B\in\B$
$$
\nu_v(B)\Eq \inf\{F(u,B):u\in\WP\}\,.\leqno(6.12)
$$
Let $B\in\B$ and $u\in\WP$ with $F(u,B)<+\infty$. In view of the
definition
of $F$ on Borel sets, by (6.11) we have
$$
F(u,B)\Eq\int_B{f_v(x,u(x))\,d\mu_v} + \nu_v(B)\Ge \nu_v(B)\,;
$$
hence, $\inf\{F(u,B):u\in\WP, F(u,B)<+\infty\}\geq \nu_v(B)$. By the
definition of $\nu_v(B)$, this implies (6.12).
\endproof
\bigskip
The following proposition shows that, given the measures $\mu$ and
$\nu$, the
function $f$ obtained in the integral representation theorem is
essentially
unique.
 
\prop{6.2}{Let $F\in\F$ with $\dF\Ne\emptyset$ and let $K=K_\om$ be
the closed
valued multifunction {}from $\om$ to $\Rm$ given by Proposition 5.1
for $A=\om$.
Let $\mu$ and $\nu$ be two positive finite Borel measures on $\om$,
with $\mu$
absolutely continuous with respect to capacity, and let
$f_1$, $f_2\colon\om\times\Rm\rightarrow [0,+\infty]$ be two Borel
functions
such that $f_1(x,\cdot)$ and $f_2(x,\cdot)$ are convex and lower
semicontinuous on $\Rm$ for $\mu\!$-a.e.\ $x\in\om$. Assume that for
every
$A\in\A$ and for every $u\in\d F(\cdot,A)$ we have $F(u,A)=\int_A
f_i(x,u(x))\,d\mu + \nu(A)$ for $i=1,2$. Then $f_1(x,\xi)=f_2(x,\xi)$ for
$\mu\!$-a.e.\ $x\in\om$ and for every $\xi\in K(x)$.}
 
\proof By a translation we can easily reduce the problem to the case
$F(0,\om)<+\infty$. Let $(u_i)$ and $C_k(x)$ be as in the proof of
Theorem~5.4.
By Proposition 2.9 we have
$$
F(u,A) = \int_A f_1(x,u(x))\,d\mu + \nu(A) = \int_A f_2(x,u(x))\,d\mu +
\nu(A)
$$
for every $A\in\A$ and for every $u\in\WP$ with $u(x)\in C_k(x)$ for
q.e.\
$x\in\om$. By Proposition 3.8 it turns out that $f_1(x,\xi)=f_2(x,\xi)$ for
$\mu\!$-a.e.\ $x\in\om$ and for every $\xi\in C_k(x)$. Hence the
equality
holds for every $\xi\in{\rm ri}K(x)\subseteq\bigcup_{k=1}^\infty
C_k(x)$,
and, therefore, for every $\xi\in K(x)$ by the continuity along line
segments
(see [28], Corollary 7.5.1).
\endproof
 
\prop{6.3}{Let $K(x)$ be a closed and convex valued multifunction
{}from
$\om$ to $\Rm$ for which
there exists a sequence $(u_k)$ of functions in $\WP$ such that
$K(x)={\rm cl}\{u_k(x):k\in\n\}$ for q.e.\ $x\in\om$. Then, there exists
a
positive finite Borel measure $\rho$ on $\om$, absolutely continuous
with
respect to capacity, such that for every $A\in\A$ and for every
$u\in\WP$
the following conditions are
equivalent:
\smallskip
\item{(i)} $u(x)\in K(x)$ for q.e.\ $x\in A$,
\smallskip
\item{(ii)} $u(x)\in K(x)$ for $\rho\!$-a.e.\ $x\in A$.
\smallskip}
 
\proof It is not restrictive to assume that $0\in K(x)$ for q.e.\
$x\in\om$.
Moreover, we can suppose that $u_k\in\WPL$ for every $k\in\n$.
Indeed, if
$T_h\circ u_k$, with $h\in\n$, denotes the truncation introduced in
Lemma~2.5,
it turns out that $K(x)={\rm cl}\{(T_h\circ u_k)(x): h, k\in\n\}$, since
$K(x)$ is
closed and convex and $(T_h\circ u_k)_h$ converges to $u_k$ strongly
in $\WP$,
as $h$ tends to $\infty$.
\par
Let us note that, by a standard cut-off argument, it is enough to consider
the case $A=\om$. Moreover, (i) clearly implies (ii) as $\rho$ is
absolutely continuous with respect to capacity.
\par
\smallskip
\noindent
{\it Step 1.} Here we prove that (ii) implies (i) for $u\in\W0$
under the additional assumption that $\partial\om$ is smooth.
\par
Let us define the convex sets
$$
\displaylines{{\cal K}=\{u\in\W0:u(x)\in K(x)\hbox{ for q.e.\ }
x\in\om\}\,,\cr
{\cal K}_k\Eq\{u\in\W0:u(x)\in K(x) +{\textstyle {1\over
k}}B_1(0)\hbox{ for
q.e.\ }x\in\om\}\,,}
$$
for every $k\in\n$. Since $\W0$ is separable, the set ${\cal K}_k$ is the
intersection of a
countable family of closed half-spaces of $\W0$. Hence, there exists a
sequence $(\mu_{k,h})_h$ in $\WPP$, with $p^\p={\textstyle {p\over
{p-1}}}$, and
a sequence $(a_{k,h})_h$ in $\r$ such that
$$
{\cal
K}_k\Eq\bigcap_{h\in\n}\{u\in\W0:\,
\langle\mu_{k,h},u\rangle\Ge a_{k,h}\,\}\,,
$$
where $\langle\cdot,\cdot\rangle$ denotes the duality pairing between
$\WPP$ and $\W0$.
\par
Denote by
$\M (\om,\Rm)$ the space of all $\Rm\!$-valued Radon measures on
$\om$ with
bounded total variation. We say that an element $T\in\WPP$ belongs to
$\M (\om,\Rm)$ if there exists $\mu\in \M (\om,\Rm)$ such that
$$
\langle T,\f\rangle\Eq \into \f \,d\mu
$$
for every $\f\in C^\infty_0(\om,\Rm)$.
In this case $T$ and $\mu$ will be identified.
\par
Let us prove that
$\mu_{k,h}\in\M (\om,\Rm)$ for every $h$,$k\in\n$. Fix $\f\in
C^\infty_0(\om,\Rm)$ with $\|\f\|_\infty\Le 1$. Since $u_1+{\textstyle
{1\over
k}}\f$ and  $u_1-{\textstyle {1\over k}}\f$ belong to ${\cal K}_k$, we
have
$$
-k(\langle\mu_{k,h},u_1\rangle - a_{k,h})\Le
\langle\mu_{k,h},\f\rangle\Le k(\langle\mu_{k,h},u_1\rangle -
a_{k,h})\,.
$$
Therefore, there exists $C_{k,h}> 0$ such that
$|\langle\mu_{k,h},\f\rangle|\leq C_{k,h}\|\f\|_\infty$ for every $\f\in
C^\infty_0(\om,\Rm)$. Hence $\mu_{k,h}|_{C^\infty_0(\om,\Rm)}$ can
be uniquely
extended to a continuous linear functional on the space of continuous
functions
on $\om$ vanishing on $\partial\om$. We conclude by the Riesz
representation
theorem.
\par
Since ${\cal K}=\bigcap_{k\in\n}{\cal K}_k$, we can assert that there
exists a
sequence $(\mu_h)$ in \break$\WPP\cap \M (\om,\Rm)$ and a
sequence $(a_h)$ in
$\r$ such that
$$
{\cal K}\Eq \bigcap_{h\in\n}\{u\in\W0: \langle\mu_h,u\rangle\Ge
a_h\}\,.\leqno(6.13)
$$
Moreover, as $\mu_h\in\WPP\cap \M (\om,\Rm)$, by [22] and [8],
Lemma~2 we
have that $|\mu_h|$ is absolutely continuous with respect to capacity.
For
every $B\in\B$ define
$$
\rho(B)\Eq \sum_{h=1}^\infty 2^{-h}{{|\mu_h|(B)}\over {|\mu_h|(\om)}}
$$
(clearly, we can assume that $|\mu_h|(\om)> 0$ for every $h\in\n$).
Then
$\rho$ is a positive finite Borel measure on $\om$ absolutely continuous
with respect to capacity. Let $g_h$ be the Radon-Nikodym derivative of
$\mu_h$ with respect to $\rho$. Then $g_h\in L^1(\om,\rho)$.
By Corollary~6 in [8], for every
$u\in\W0\cap\LL$, we have $u\cdot g_h\in L^1(\om,\rho)$ and
$$
\langle\mu_h,u\rangle\Eq \into u\cdot g_h\,d\rho\,.\leqno(6.14)
$$
Let us now prove that (ii) implies (i).
Fix $u\in\W0$ with $u(x)\in K(x)$ for $\rho$-a.e.\ $x\in\om$. Assume
first
that $u$ belongs to $\LL$.
\par
For every $h\in\n$ and $v\in{\cal K}\cap\LL$, by (6.14) we have
$$
a_h\Le \langle\mu_h,v\rangle\Eq\into v\cdot g_h\,d\rho\,,
$$
hence
$$
a_h\Le \inf_{v\in {\cal K}\cap L^\infty}\into v\cdot
g_h\,d\rho\,.\leqno(6.15)
$$
In view of the fact that the functions $u_k$ are in $\LL$, it turns out
that
$K(x)=\rho\hyphen{\displaystyle\essup_{v\in
{\cal K}\cap L^\infty}}\,\{v(x)\}$ for $\rho$-a.e.\ $x\in\om$; then
Theorem~4.4
yields
$$
\inf_{v\in {\cal K}\cap L^\infty}\into v\cdot g_h\,d\rho
\Eq\into\inf_{\xi\in
K(x)}\xi\cdot g_h(x)\,d\rho\,.\leqno(6.16) $$
By the assumption $u(x)\in K(x)$ for $\rho$-a.e.\ $x\in\om$,
by (6.15) and (6.16), it follows that
$$
a_h\Le \into u(x)\cdot g_h(x)\,d\rho\Eq \langle\mu_h,u\rangle
$$
for every $h\in\n$. By (6.13) this proves that $u\in {\cal K}$, i.e.,
$u(x)\in
K(x)$ for q.e.\ $x\in\om$.
\par
Consider now a general $u\in\W0$; let us note that, since $K(x)$
is convex and $0\in K(x)$ for q.e.\  $x\in\om$, the condition $u(x)\in
K(x)$
for $\rho$-a.e.\ $x\in\om$ implies that
$(T_h\circ u)(x)\in K(x)$ for $\rho$-a.e.\ $x\in\om$ and for every
$h\in\n$.
The previous step and the q.e.\ convergence of $(T_h\circ u)$ to $u$
allow us to
conclude as $K(x)$ is closed.
\par
\smallskip
\noindent
{\it Step 2.} Let us now prove that (ii) implies (i) for every
$u\in\WP$ without assuming the smoothness of the boundary of $\om$.
\par
Let $(\om_h)$ be a sequence of open subsets of $\om$ with
$\om_h\subset\subset\om_{h+1}$, $\bigcup_h \om_h=\om$, and
$\partial\om_h$
smooth. Let $\f_h$ be a $C^1_0(\om_h)$ function with $\f_h=1$ on
$\om_{h-1}$ and $0\leq\f_h\leq 1$. Define
$$
K_h(x)\Eq {\rm cl}\{\f_h(x)u_k(x): k\in\n\}
$$
for q.e.\ $x\in\om_h$. By Step~1 there exists a positive finite Borel
measure $\rho_h$ on $\om_h$, absolutely continuous with respect to
capacity, and such that for every $u\in\W0$ the condition $u(x)\in
K_h(x)$ for
q.e.\ $x\in\om_h$ is equivalent to the condition $u(x)\in K_h(x)$ for
$\rho_h$-a.e.\ $x\in\om_h$. We can consider $\rho_h$ as a measure on
$\om$ by
setting $\rho_h(B)=\rho_h(B\cap\om_h)$ for $B\in\B$. Let us define
$$
\rho(B)\Eq \sum_{h=1}^\infty 2^{-h}{{\rho_h(B)}\over{\rho_h(\om)}}
$$
for every $B\in\B$. Then $\rho$ is a positive finite Borel measure on
$\om$
which is absolutely continuous with respect to capacity.
\par
Let us fix $u\in\WP$ with
$u(x)\in K(x)$ for $\rho$-a.e.\ $x\in\om$. Then, $\f_h(x)u(x)\in K_h(x)$
for
$\rho_h$-a.e.\ $x\in\om_h$, so that $\f_h(x)u(x)\in K_h(x)$ for
q.e.\ $x\in\om_h$. Since $\f_h=1$ on $\om_{h-1}$, we obtain that
$u(x)\in K(x)$
for q.e.\ $x\in\om_{h-1}$. As $h$ is arbitrary, we conclude that
$u(x)\in K(x)$ for q.e.\ $x\in\om$.
\endproof
 
\lemma{6.4}{Let $F\in\F$ and define
$\om_0$ to be the union of all $A\in\A$ such that
$\d F(\cdot,A)\ne\emptyset$.
Then $\d F(\cdot,A)\ne\emptyset$ for every $A\in\A$ with
$A\subset\subset\om_0$.}
 
\proof By induction we can reduce ourselves to prove that, whenever
$A_1$, $A_2$ are open subsets of $\om$ with $\d
F(\cdot,A_1)\ne\emptyset$
and
$\d F(\cdot,A_2)\ne\emptyset$, then $\d F(\cdot,A)\ne\emptyset$ for
every open set $A\subset\subset A_1\cup A_2$.
\par
Let $A_1$, $A_2$ and $A$ be as above, and let $A^\p\subset\subset
A_1$ with
$A\subset\subset A^\p\cup A_2$. We shall show that $\d
F(\cdot,A^\p\cup
A_2)\ne\emptyset$, which clearly implies $\d
F(\cdot,A)\ne\emptyset$.
Consider a function $\f\in C^1_0(A_1)$ with $\f=1$ on $A^\p$ and
$0\le\f\le 1$. By assumption we can find $u\in \d F(\cdot,A_1)$ and
$v\in\d F(\cdot,A_2)$; define $w= \f u+(1-\f)v$. Then, by the usual
properties of the class $\F$ it is easy to see that $F(w,A^\p\cup
A_2)<+\infty$.
\endproof
 
\th{6.5}{Let $\om$ be an open subset of $\r^n$ (not necessarily
bounded) and let $F\in\F$. Then,
there exist a positive finite Borel measure $\mu$ on $\om$, absolutely
continuous with respect to capacity, a positive Borel measure $\nu$ on
$\om$,
and a Borel function $f\colon\om\times\Rm\rightarrow[0,+\infty]$
with the following properties:
\smallskip
\item{(i)} for every $x\in\om$ the function $f(x,\cdot)$ is convex and
lower semicontinuous on $\Rm$;
\smallskip
\item{(ii)} for every $u\in\WP$ and for every $A\in\A$
\smallskip\noindent
$$
F(u,A)\Eq
\int_A{f(x,u(x))\,d\mu} + \nu(A)\,.\leqno(6.17)
$$
}
 
\proof Let $\overline\nu$ be the function on $\B$ defined in (6.1).
Let us set
$$
\nu(A)\Eq \Sup\{\overline\nu(A^\p): A^\p\in\A,\ A^\p\subset\subset
A\}
$$
for every $A\in\A$. Clearly $\nu$ is increasing with respect to
the inclusion, and $\nu(\emptyset)=0$. Moreover, for every $A\in\A$
$$
\nu(A)\L +\infty\ \Rightarrow\ A\subseteq\om_0\,.\leqno(6.18)
$$
Clearly $\nu$ is inner regular on $\A$; therefore, by
Proposition 5.5 and Theorem~5.6 in~[19], to prove that $\nu$ can
be extended to a Borel measure on $\om$ it suffices to show that
$\nu$ is subadditive and superadditive on $\A$.
\par
Let $A_1$, $A_2\in\A$, and note that
$$
\nu(A_1\cup A_2)\Eq \Sup\{\overline\nu(A_1^\p\cup A_2^\p):
A_i^\p\in\A,\ A_i^\p\subset\subset A_i\ \ (i=1,2)\}\,.\leqno(6.19)
$$
If $\nu(A_1)$ and $\nu(A_2)$ are finite, then
$A_1,A_2\subseteq\om_0$
by (6.18). Since, by Lemma~6.4, $\overline\nu$ is a measure on every
$\om^\p\subset\subset\om_0$, {}from (6.19) it follows that
$\nu(A_1\cup A_2)
\le \nu(A_1)+\nu(A_2)$. In a similar way we get superadditivity.
\par
This allows us to conclude that the set function
$\nu\colon\B\rightarrow[0,+\infty]$ defined by
$$
\nu(B)\Eq\inf\{\nu(A):A\in\A,\ B\subseteq A\}
$$
is a Borel measure on $\om$ and that $\nu (B)={\overline\nu}(B)$ for
every
$B\in\B$ with $B\subset\subset\om_0$.
\bigskip
Let us construct now $\mu$ and $f$. Let $(\om_h)$ be a sequence
of open subsets of $\om_0$ with a smooth boundary such that
$\om_h\subset\subset\om_{h+1}$ and $\om_0=\cup_h \om_h$.
In particular, $\d F(\cdot,\om_h)\ne\emptyset$ by Lemma~6.4.
Since there
exists an extension operator {}from $W^{1,p}(\om_h,\Rm)$
to $\WP$, it is possible to apply Theorem~6.1 to each $\om_h$ using
$\WP$ instead of $W^{1,p}(\om_h,\Rm)$.
Let $K_{\om_h}$ be the multifunction defined in Proposition 5.1
for $A=\om_h$;
then there exist a
positive finite Borel measure $\mu_h$ on $\om_h$, absolutely
continuous with
respect to capacity, and a Borel function $f_h\colon
\om_h\times\Rm\rightarrow[0,+\infty]$ such that
\smallskip
\item{(a)} for every $x\in \om_h$ the function $f_h(x,\cdot)$ is convex
and
lower semicontinuous on $\Rm$;
\smallskip
\item{(b)} for every $u\in\WP$ and $A\in{\cal A}(\om_h)$
$$
F(u,A)\Eq
\cases{
\int_A{f_h(x,u(x))\,d\mu_h} + \overline\nu(A)\,,& if
                   \ $u(x)\in K_{\om_h}(x)$\quad\hbox{for q.e.\ $x\in A\,,$}
\cr
\Extra
\cr
+\infty\,,& otherwise$\,.$\cr}
$$
\smallskip
Moreover, we have a uniqueness property for the integrand
as stated in Proposition 6.2.
\par
By Proposition 6.3, for every $h\in\n$
there exists a positive finite Borel measure $\rho_h$ on
$\om_h$, absolutely continuous with respect to capacity, such that for
every
$u\in\WP$ the condition $u(x)\in K_{\om_h}(x)$ for q.e.\ $x\in \om_h$
is
equivalent to the condition $u(x)\in K_{\om_h}(x)$ for $\rho_h\!$-a.e.\
$x\in
\om_h$. Let us consider $\mu_h$ and $\rho_h$ as measures on $\om$
by setting $\mu_h(B)=\mu_h(B\cap \om_h)$ and
$\rho_h(B)=\rho_h(B\cap \om_h)$
for every $B\in\B$. Define
$$
\mu(B)\Eq\sum_{h=1}^\infty 2^{-h}{{(\mu_h+\rho_h)(B)}\over{(\mu_h +
\rho_h)(\om)}}\,, $$
$$
g_h(x,\xi)\Eq
\cases{
f_h(x,\xi)\,{{d\mu_h}\over {d\mu}}(x)\,,
& if $x\in\om_h$ and $\xi\in K_{\om_h}(x)\,,$\cr
\Extra
\cr
+\infty\,,& otherwise,\cr}
$$
where $d\mu_h/d\mu$ is a fixed Borel representative of
the Radon-Nikodym derivative. Then $\mu$ is a positive
finite Borel measure on $\om$, absolutely continuous with respect to
capacity.
Since for every $h\in\n$ there is a sequence $(v_i)$ in $\WP$ such
that $K_{\om_h}(x)={\rm cl}\{v_i(x):i\in\n\}$, by Theorem~III.9 and
Proposition~III.13 in~[10], the graph of $K_{\om_h}$
belongs to ${\cal B}(\om_h)\times{\cal B}(\Rm)$. Therefore
$g_h\colon\om_h\times\Rm\rightarrow[0,+\infty]$
is a Borel function, and for every $x\in \om_h$, the function
$g_h(x,\cdot)$ is convex and lower
semicontinuous on $\Rm$.
\par
By recalling that $K_{\om_h}=K_{\om_{h+1}}$ q.e.\ on $\om_h$ (see
Remark~5.2), and by using the uniqueness property of the integrand
mentioned above, we easily obtain that
$g_h(x,\cdot)=g_{h+1}(x,\cdot)$
for $\mu\!$-a.e.\ $x\in\om_h$.
\par
Therefore, there exists a Borel function $f\colon\om\times\Rm
\rightarrow [0,+\infty]$ satifying (i) and such that for every $h\in\n$
$$
f(x,\cdot)=g_h(x,\cdot)\qquad\hbox{on $\Rm$ for $\mu\!$-a.e.\
$x\in\om_h\,.$}
$$
\par
Let us now prove (ii). Fix $u\in\WP$ and $A\in\A$. If $A\setminus
\om_0\ne\emptyset$, then $\nu(A)=+\infty$ by (6.18). On the other
hand, by
the definition of $\om_0$ we have $F(u,A)=+\infty$  for every
$u\in\WP$.
Therefore
$$
F(u,A)\Eq\int_A f(x,u(x))\,d\mu + \nu(A)\,.
$$
\par
Let now $A\subseteq\om_0$ and $A^\p\subset\subset A$. Then
$A^\p\subseteq\om_h$ for a suitable $h\in\n$. In view of the
properties
of the measure $\rho_h$, {}from the definition of $g_h$ and $f$ we
easily
obtain
$$
F(u,A^\p)\Eq\int_{A^\p} g_h(x,u(x))\,d\mu + \overline\nu(A^\p)
\Eq\int_{A^\p} f(x,u(x))\,d\mu + \overline\nu(A^\p)\,.
$$
Therefore (ii) follows {}from the definition of $\nu$ taking the
supremum
for $A^\p\subset\subset A$.
\endproof
 
\rem{6.6}{Let $F\in\F$ and $\om_0$ be as in Lemma~6.4.
By Proposition~5.1 and Remark~5.2,
there exists a closed valued multifunction $K$ {}from $\om_0$ to $\Rm$,
unique up to sets of capacity zero, such that
$$
K(x)\Eq K_A(x)\qquad\hbox{for q.e.\ $x\in A$}\leqno(6.20)
$$
whenever $A\in\A$ and $\d F(\cdot,A)\ne\emptyset$.
Moreover, $K(x)$ is
non-empty and convex for q.e.\ $x\in\om_0$.
 
	It is clear that the function $f$ constructed in the proof of
Theorem~6.5 satisfies the additional condition
$f(x,\xi)=+\infty$ for $\mu$-a.e.\ $x\in\Omega_0$ and for every
$\xi\notin K(x)$. This is not necessarily true for every function $f$
which satisfies conditions (i) and (ii) of Theorem~6.5.
	Let us consider, for instance, the functional
$$
F(u,A)=\cases{0\,,&if $u=0$ a.e.\ on $A$,\cr
+\infty\,,&otherwise,\cr}
$$
in the case $n=m=1$ and $\Omega=\r$. Then, clearly, $K(x)=\{0\}$ for
q.e.\ $x\in\r$ and
$$
F(u,A)=\int_A f(x,u(x))\,dx\,,\leqno(6.21)
$$
with
$$
f(x,\xi)=\cases{0\,,&if $\xi=0$,\cr
+\infty\,,&if $\xi\neq0$.\cr}
$$
But (6.21) holds also with $f(x,\xi)=a(x)|\xi|^2$, where
$a\colon\r\to{[0,+\infty[}$ is any finite valued Borel function such that
$\int_A a(x)\,dx=+\infty$ for every open subset $A$ of~$\r$ (see~[23],
Section~43, Exercise~7).
}
 
\rem{6.7}{If $\dF=\emptyset$ and $\nu$ is not necessarily finite, the
uniqueness result of Proposition~6.2
still holds,
with an obvious localization of the proof, in the weaker
form:
$$f_1(x,\xi)=f_2(x,\xi)\quad\hbox{for $\mu$-a.e.\ $x\in\om_0$ and
for every
$\xi\in K(x)\,,$}
$$
where $\om_0$ is defined in Lemma~6.4 and $K(x)$ is now defined
by~(6.20).}
\bigskip
 
\parag{7}{Quadratic functionals}
 
In this section we show how certain algebraic properties of the
functional $F$ are inherited by the integrand which appears in
the representation of $F$ according to Theorem~6.5. We recall that a
cone in a vector space $X$ (with vertex at~$0$) is a set $K$ such that
$tx\in K$ for every $t>0$ and for every $x\in K$.
 
\defin{7.1}{Let $X$ be a real vector space and let $p\in \r$.
We say that a function $f\colon X\rightarrow [0,+\infty]$ is:
\smallskip
\item{(i)} {\it positively homogeneous of degree\/} $p$ on a cone $K$ if
$F(tx)=t^pF(x)$ for every $t>0$ and for every $x\in K$;
\smallskip
\item{(ii)} {a (non-negative) {\it quadratic form}
(with extended real values) on $X$ if there exist a linear subspace $Y$ of
$X$ and a symmetric bilinear form $B\colon Y\times Y\rightarrow \r$
such that
$$
F(x)\Eq\cases{B(x,x)\,,& if $x\in Y$,\cr
+\infty\,,& if $x\in X\setminus Y$.\cr}
$$
}
\vskip-12truept
\item{} We shall refer to $Y$ as the {\it domain} of $F$.
}
 
\rem{7.2}{In the previous definition it is not restrictive to assume
that $B$ is
defined over all of $X\times X$. Indeed, let $Z$ be an
algebraic complement of $Y$ in $X$ and denote by $P\colon
X\rightarrow Y$ the
canonical projection on $Y$ associated to the pair $(Y,Z)$. Then, it is
enough to consider the extension $(x,y)\mapsto B(Px,Py)$ defined for
every
$(x,y)\in X\times X$. As a consequence, if $X$ is finite dimensional and
${\rm dim}X=m$, then there exists an $m\times m$ symmetric matrix
$(a_{ij})$
such that $F(x)=\sum_{i,j=1}^m a_{ij}x^ix^j$ for every $x\in Y$, where
$x^1,\ldots, x^m$ denote the components of $x$ with respect to
a fixed basis of $X$.
}
 
\th{7.3}{Let $F\in\F$ and $p\in\r$. Let $f$, $\mu$, $\nu$ be as in
Theorem~6.5 and let $\om_0$ and $K$ be as in Remark~6.6. Assume
that $f(x,\xi)=+\infty$ for $\mu$-a.e.\ $x\in\om_0$ and for every
$\xi\notin K(x)$. Then the following properties hold:
\smallskip
\item{(i)} if $F(\cdot,A)$ is
positively homogeneous of degree $p$ on $\WP$ for every $A\in\A$,
then $K(x)$
is a closed convex cone for q.e.\ $x\in\om_0$, and $f(x,\cdot)$ is
positively
homogeneous of degree $p$ on $K(x)$ for $\mu$-a.e.\ $x\in\om_0$; if,
in
addition, $p\ne 0$, then $\nu (B)=0$ for every $B\in{\cal B}(\om_0)$;
\smallskip
\item{(ii)} if $F(\cdot,A)$ is a quadratic
form on $\WP$ for every $A\in\A$,
then $\nu=0$,
$K(x)$ is a linear subspace of $\Rm$ for q.e.\ $x\in\om$, and for
$\mu$-a.e.\ $x\in\om$ the
function $f(x,\cdot)$ is a quadratic form on $\Rm$ with domain $K(x)$.
\smallskip}
 
\proof Proof of (i).\ \ For every $A\in\A$ the positive homogeneity of
degree $p$ implies that $tu\in
\d F(\cdot,A)$ whenever $t>0$ and $u\in \d F(\cdot,A)$. Recalling the
definition and properties of $K_A$ given in Proposition~5.1, it is easy to
see
that, if ${\rm dom}F(\cdot,A)\Ne\emptyset$, then there exists a set
$N\subseteq\om$ with ${\rm cap}(N)=0$ and such that for every $x\in
A\setminus
N,\,q\in\q^+\setminus\{0\}$ and $\xi\in K_A(x)$ we have $q\xi\in
K_A(x)$. Since
$K_A(x)$ is closed, it follows that $K_A(x)$ is a cone for every $x\in
A\setminus N$. The convexity of $K_A(x)$ is proved in Proposition
5.1(iii). By
the definition (6.20) of $K(x)$ and by the definition of $\om_0$ we
conclude
that $K(x)$ is a closed convex cone for q.e.\ $x\in\om_0$.
\par
Let us now prove that
$f(x,\cdot)$ is positively homogeneous of degree $p$ on $K(x)$ for
$\mu$-a.e.\
$x\in\om_0$. Let us first consider the case $p=0$. If $A\in\A$ and
$u\in
\d F(\cdot,A)$, then the function $t\mapsto F(tu,A)$ {}from $[0,1]$ into
$[0,+\infty]$ is convex and lower semicontinuous; moreover,
$F(tu,A)=F(u,A)<+\infty$ for every $t>0$. Therefore,
$F(0,A)=\displaystyle\lim_{t\to 0^+} F(tu,A)=F(u,A)$ for every $A\in\A$
and $u\in\d F(\cdot,A)$. This shows that $0\in K(x)$ for q.e.\
$x\in\om_0$. By the uniqueness of the integrand
stated in Remark~6.7, we conclude that $f(x,\xi)=f(x,0)$ for
$\mu$-a.e.\ $x\in\om_0$ and
for every $\xi\in K(x)$.
\par
Assume now $p\ne 0$. Let $A\in\A$ and $w\in \d F(\cdot,A)$ . Since
$\nu(A)\leq
F(tw,A)=t^pF(w,A)$ for every $t>0$, taking the limit as $t\to 0^+$ or $t\to
+\infty$ according to whether $p>0$ or $p<0$, we get $\nu(A)=0$. In
view of
the positive homogeneity of $F(\cdot,A)$ and (6.17) we have
$$
F(u,A)\Eq\int_A {1\over {t^p}}f(x,tu(x))\,d\mu
$$
for every $A\in\A$, $u\in\WP$ and $t>0$. Now, by the uniqueness of
the
integrand (Remark 6.7), we have $f(x,\xi)=(1/t^p)f(x,t\xi)$ for
$\mu$-a.e.\ $x\in\om_0$ and for every $\xi\in K(x)$.
\medskip
\par\noindent
Proof of (ii).\ \ Assume that $F(\cdot,A)$ is a quadratic form for every
$A\in\A$. Then $F(0,A)=0$ for every $A\in\A$, hence $\om_0=\om$,
$\nu(B)=0$
for every $B\in\B$, and $f(x,0)=0$ for $\mu\!$-a.e.\ $x\in\om$. Directly
{}from Definition 7.1 it follows that $\d F(\cdot,A)$ is a linear space; in
particular, $u+v$ and $-u$ belong to $\d F(\cdot,A)$ whenever
$u,v\in\d
F(\cdot,A)$. As in the first part of (i), it can be shown that for
q.e.\ $x\in\om$, $\xi+\eta$ and $-\xi$ belong to $K(x)$ if $\xi, \eta\in
K(x)$. Since $K(x)$ is a cone (part (i)), this guarantees that
$K(x)$ is a linear subspace of $\Rm$ for q.e.\ $x\in\om$.
\par
If $X$ is a
(real) vector space, it is well known (Fr\'echet-Von~Neumann-Jordan
Theorem, see, for instance, [30]) that a function $F\colon
X\rightarrow
[0,+\infty]$ is a quadratic form if and only if $F(0)=0$, $F$ is positively
homogeneous of degree $2$, and satisfies the following ``parallelogram
identity":
$$
F(\xi+\eta)+F(\xi-\eta)\Eq 2F(\xi) + 2F(\eta)
$$
for every $\xi,\eta\in X$. Since $f(x,0)=0$ for $\mu\!$-a.e.\ $x\in\om$
and
$f(x,t\xi)=t^2f(x,\xi)$ for $\mu\!$-a.e.\ $x\in\om$ and for every $t>0$,
$\xi\in K(x)$ (see part (i)), to complete the proof of (ii) it remains only
to show that $f(x,\cdot)$ satisfies the parallelogram identity on $K(x)$
for
$\mu$-a.e.\ $x\in\om$. Define the functional $G\colon
[\WP]^2\times\A\rightarrow [0,+\infty]$ as
$$
G(u,v,A)=F(u+v,A)+F(u-v,A)=2F(u,A)+2F(v,A)\,.
$$
Since $\nu =0$, {}from (6.17) we
obtain
$$
\displaylines{\llap(7.1)\hfill G(u,v,A)\Eq\int_A
2[f(x,u(x))+f(x,v(x))]\,d\mu
\Eq\hfill
\cr\int_A[f(x,u(x)+v(x))+f(x,u(x)-v(x))]\,d\mu}
$$
for every $A\in\A$. Since $[\WP]^2$ can be identified with
$W^{1,p}(\om,\r^{2m})$, we can apply Remark~6.7 to the functional $G$,
with
the set $K(x)\times K(x)$ playing the role of $K(x)$ for q.e.\ $x\in\om$.
Therefore, (7.1) gives that
$$
2(f(x,\xi)+f(x,\eta))=f(x,\xi+\eta)+f(x,\xi-\eta)
$$
for $\mu$-a.e.\ $x\in\om_0$ and for every $(\xi,\eta)\in K(x)\times
K(x)$.
\endproof
 
\cor{7.4}{
Let $F\in\F$. Assume that $F(\cdot,A)$ is a quadratic form on
$\WP$ for every $A\in\A$. Then there exist:
\smallskip
\item {(i)} a positive finite
Borel measure $\mu$ on $\om$, absolutely continuous with respect to
capacity,
\smallskip
\item {(ii)} a symmetric $m\times m$ matrix $(a_{ij})$ of Borel
functions
{}from $\om$ into $\r$ such that $\sum_{i,j=1}^m a_{ij}(x)\xi^i\xi^j\ge
0$ for q.e.\ $x\in\om$ and for every $\xi\in\Rm$,
\smallskip
\item {(iii)} for every $x\in\om$ a linear subspace
$V(x)$ of $\Rm$,
\smallskip
\noindent
with the following properties: for every $u\in\WP$ and $A\in{\cal
A}(\om)$
\smallskip
\item {(a)} if $F(u,A)<+\infty$, then $u(x)\in V(x)$ for q.e.\ $x\in A$;
\smallskip
\item {(b)} if $u(x)\in V(x)$ for q.e.\ $x\in A$, then
$F(u,A)\Eq\int_A \sum_{i,j=1}^m a_{ij}(x)u^i(x)u^j(x)\,d\mu$.
\smallskip
}
 
\proof The conclusion follows {}from
Theorems~6.5 and 7.3(ii), and from Remark~7.2.
\endproof

\goodbreak\bigskip\bigskip\goodbreak
\centerline{\bf References}
\nobreak\bigskip\nobreak
 
\baselineskip=12pt
{\eightpoint
\frenchspacing
\smallskip\item{[1]} ALBERTI G.: Integral respresentation of local
functionals. Preprint Scuola Normale, Pisa, 1990.
 
\smallskip\item{[2]} ALBERTI G., BUTTAZZO G.: Integral respresentation
of functionals defined on Sobolev spaces. {\it
Composite Media and Homogenization Theory (Trieste, 1990)\/}, 1-12,
{\it Birkh\"auser, Boston\/}, 1991.
 
\smallskip\item{[3]} AMBROSIO L., DAL MASO G.: On the relaxation in
$BV(\Omega;{\bf R}^m)$ of quasi-convex integrals.
{\it J. Funct. Anal.\/}, to appear.
 
\smallskip\item{[4]} AMBROSIO L., PALLARA D.: Integral
representation of relaxed functionals on $BV({\bf R}^n;{\bf R}^k)$ and
polyhedral  approximation.
Preprint Scuola Normale, Pisa, submitted to
{\it Indiana Univ. Math.~J.\/}
 
\smallskip\item{[5]} ATTOUCH H., PICARD C.: Variational
inequalities with varying obstacles: the general form of the limit
problem. {\it J. Funct. Anal.} {\bf 50} (1983), 329-386.
 
\smallskip\item{[6]}  BOUCHITT\'E G., DAL MASO G.: Integral
representation and relaxation of convex local functionals on
BV($\Omega$).
Preprint SISSA, Trieste, 1991.
 
\smallskip\item{[7]}  BOUCHITT\'E G. , VALADIER M.: Integral
representation of convex functionals on a space of measures. {\it
J. Funct. Anal.} {\bf 80} (1988), 398-420.
 
\smallskip\item{[8]} BREZIS H., BROWDER F.E.: Some properties of
higher
order Sobolev spaces. {\it J. Math. Pures Appl. (9)} {\bf 61} (1982),
245-259.
 
\smallskip\item{[9]}  BUTTAZZO G.: Semicontinuity, Relaxation and
Integral Representation Problems in the Calculus of Variations. Pitman
Res. Notes in Math. 207, Longman, Harlow, 1989.
 
\smallskip\item{[10]}  CASTAING C., VALADIER M.: Convex Analysis
and
Measurable Multifunctions. Lecture Notes in Math. 580, Springer-Verlag,
Berlin, 1977.
 
\smallskip\item{[11]}  CHOQUET G.: Lectures on Analysis. Vol. I.
Mathematics Lecture Note Series, Benjamin, New York, 1969.
 
\smallskip\item{[12]} DAL MASO G.: Asymptotic behaviour of minimum
problems with bilateral obstacles. {\it Ann. Mat. Pura Appl. (4)}
{\bf 129} (1981), 327-366.
 
\smallskip\item{[13]} DAL MASO G.: On the integral representation
of certain local functionals.
{\it Ricerche Mat.} {\bf 32} (1983), 85-113.
 
\smallskip\item{[14]} DAL MASO G., DEFRANCESCHI A., VITALI E.:
A characterization of $C^1\!$-convex sets in Sobolev spaces.
{\it Manuscripta Math.\/}, to appear.
 
\smallskip\item{[15]} DAL MASO G., LONGO P.: $\Gamma$-limits
of obstacles. {\it Ann. Mat. Pura Appl. (4)}  {\bf 128}
(1980), 1-50.
 
\smallskip\item{[16]} DAL MASO G., PADERNI G.: Integral
representation
of some convex local functionals.
{\it Ricerche Mat.} {\bf 36} (1987), 197-214.
 
\smallskip\item{[17]} DEFRANCESCHI A., VITALI E.: Paper in
preparation.
 
\smallskip\item{[18]} DE GIORGI E., DAL MASO G., LONGO P.:
$\Gamma$-limiti di ostacoli. {\it Atti Accad. Naz. Lincei Rend. Cl. Sci.
Fis. Mat. Natur. (8)} {\bf 68} (1980), 481-487.
 
\smallskip\item{[19]} DE GIORGI E., LETTA G.: Une notion g\'en\'erale
de
convergence faible pour des fonctions croissantes d'ensemble. {\it Ann.
Scuola Norm. Sup. Pisa Cl. Sci. (4)} {\bf 4} (1977), 61-99.
 
\smallskip\item{[20]}  FEDERER H., ZIEMER W.P.: The Lebesgue set
of a function whose distribution derivatives are {\it p}th power
summable. {\it Indiana Univ. Math. J.} {\bf 22} (1982), 139-158.
 
\smallskip\item{[21]} FONSECA I., RYBKA P.: Relaxation of multiple
integrals
in the space  $BV(\Omega;{\bf R}^p)$. Research Report, Carnegie Mellon
University, Pittsburgh, 1991.
 
\smallskip\item{[22]} GRUN-REHOMME M.: Caract\'erisation du
sous-diff\'erentiel d'int\'egrandes convexes dans les espaces de Sobolev.
{\it
J. Math. Pures Appl.}  {\bf 56} (1977), 149-156.
 
\smallskip\item{[23]} HALMOS  P.R.: Measure Theory. Van Nostrand,
Princeton, 1950.
 
\smallskip\item{[24]} LOJASIEWICZ S., Jr.: Parametrizations of
convex sets. {\it J. Approx. Theory}, to appear.
 
\smallskip\item{[25]}  MORREY C.B., Jr.: Multiple Integrals in the
Calculus of Variations. Springer-Verlag,  Berlin, 1966.
 
\smallskip\item{[26]} ORNELAS A.: Parametrization of Carath\'eodory
multifunctions. {\it Rend. Sem. Mat. Univ. Padova}
{\bf 83} (1990), 33-44.
 
\smallskip\item{[27]} PICARD C.  Probl\`eme biharmonique avec
obstacles
variables. Th\`ese, Univ. Paris-Sud, 1984.
 
\smallskip\item{[28]}  ROCKAFELLAR R.T.: Convex Analysis. Princeton
University Press, Priceton, 1970.
 
\smallskip\item{[29]} VALADIER M.: Multi-applications mesurables \`a
valeurs convexes compactes. {\it J. Math. Pures Appl.}  {\bf 50} (1971),
265-297.
 
\smallskip\item{[30]} YOSIDA  K.: Functional Analysis.
Springer-Verlag,  Berlin, 1965.
 
\smallskip\item{[31]}  ZIEMER W.P.: Weakly Differentiable Functions.
Springer-Verlag,  Berlin, 1989.
}
 
\bigskip
\bigskip
\baselineskip=11pt
{\settabs
\+\quad\quad &Gianni DAL MASO\quad\quad
&Anneliese DEFRANCESCHI\quad\quad &Enrico VITALI
\cr
\+&Gianni DAL MASO &Anneliese DEFRANCESCHI&Enrico VITALI
\cr
\+&SISSA &Universit\`a di Trento &Universit\`a di Parma
\cr
\+&Via Beirut 2/4 &Dip.\ di Matematica & Dip.\ di Matematica
\cr
\+&I-34014 TRIESTE &I-38050 POVO &Via D'Azeglio 85/A
\cr
\+&&(TRENTO)&I-43100 PARMA
\cr}
\bye